\documentclass[10pt,twoside]{article}
\usepackage[colorlinks,citecolor=blue]{hyperref}
\usepackage[hyperref,frak]{paper_diening}
\usepackage[]{amsmath}
\usepackage{exscale}
\usepackage{amssymb, bm}
\usepackage{esint}
\usepackage[mathscr]{eucal}
\usepackage[amsmath,thmmarks]{ntheorem}
\usepackage{bbm}
\usepackage[hang,small,bf]{caption}
\usepackage{fancyhdr}
\usepackage{color}
\usepackage{abstract}
\usepackage{mathtools}
\usepackage{enumitem}

\theoremseparator{:}
\newtheorem{theorem}{Theorem}[section]

\newtheorem{lemma}[theorem]{Lemma}

\newtheorem{proposition}[theorem]{Proposition}
\newtheorem{corollary}[theorem]{Corollary}

\theorembodyfont{\normalfont}
\newtheorem{remark}[theorem]{Remark}

\newtheorem*{acknow}{Acknowledgements}
\theoremsymbol{\ensuremath{\square}}
\theoremheaderfont{\normalfont\scshape}
\theorembodyfont{\normalfont}
\newtheorem*{proof}{Proof}

\newcommand{\R}{\ensuremath{\mathbb{R}}}
\newcommand{\N}{\ensuremath{\mathbb{N}}}

\newcommand{\dd}{\mathrm{d}}

\renewcommand{\epsilon}{\varepsilon}

\usepackage[active]{srcltx}

\usepackage{titlesec}
\titleformat{\section}{\bf\large}{\thesection}{1em}{}
\titleformat{\subsection}{\bf\normalsize}{\thesubsection}{1em}{}

\usepackage[top=2.9cm,bottom=2.8cm,left=2.9cm,right=2.9cm]{geometry}
\theoremlisttype{all}
\newcommand{\comment}[1]{}
\newcommand{\longversion}[1]{}

\renewcommand{\Xi}{X}

\renewcommand{\phi}{\varphi}

\newcommand{\RN}{\setR^N}

\newcommand{\kappaA}{{\kappa_{\!\mathcal{A}}}}

\newcommand{\meantmp}[2]{#1\langle{#2}#1\rangle}

\def\la1{\lambda_1}

\newcommand{\F}{\mathcal{F}}

\newcommand{\kabs}[1]{\ensuremath{\vert#1\vert}}

\newcommand{\mint} {\mathop{\int\hskip -1,05em -\,}}
\newcommand{\mI}[1]{\mint\nolimits_{\!\!\!\!#1}}

\numberwithin{equation}{section}

\newtagform{boldbrackets}{\hspace{4ex}\bf(}{\bf)}
\newtagform{standardbrackets}{)}{)}



\makeatletter
\let\orgdescriptionlabel\descriptionlabel
\renewcommand*{\descriptionlabel}[1]{%
  \let\orglabel\label
  \let\label\@gobble
  \phantomsection
  \edef\@currentlabel{#1}%
  \let\label\orglabel
  \orgdescriptionlabel{#1}%
}
\makeatother

\allowdisplaybreaks[3] 

\color{black}

\title{\vspace{-1cm}
{{\bf Partial regularity for Minimizers of Discontinuous Quasiconvex Integrals \\ with General Growth}}
}
\date{}
\author{Christopher Goodrich, Giovanni Scilla and Bianca Stroffolini}
\newcommand{\Addresses}{{
  \bigskip
  \footnotesize
C. S.~Goodrich, \textsc{School of Mathematics and Statistics, UNSW Sydney, Sydney, NSW 2052, Australia}\par\nopagebreak
  \textit{E-mail address}, C. S.~Goodrich: \texttt{c.goodrich@unsw.edu.au}\\

G.~Scilla, \textsc{Dipartimento di Scienze di Base ed Applicate per l'Ingegneria (SBAI), Sapienza Universit\`{a} di Roma, Via A. Scarpa 16, 00169 Roma, Italy}\par\nopagebreak
  \textit{E-mail address}, G.~Scilla: \texttt{giovanni.scilla@uniroma1.it}\\

B.~Stroffolini, \textsc{Dipartimento di Ingegneria Elettrica e delle Tecnologie dell'Informazione , Universit\`{a} di Napoli Federico II, Via Claudio, 80125 Napoli, Italy}\par\nopagebreak
  \textit{E-mail address}, B.~Stroffolini: \texttt{bstroffo@unina.it}
}}

\begin{document}

\maketitle

\begin{abstract} We prove the partial H\"older continuity for minimizers of quasiconvex functionals 
\begin{equation*}
\F({\bf u}) \colon =\int_{\Omega} f(x,{\bf u},D{\bf u})\,\dd x,
\end{equation*}
where $f$ satisfies a uniform VMO condition with respect to the $x$-variable and is continuous with respect to ${\bf u}$. The growth condition with respect to the gradient variable is assumed a general one.
\end{abstract}


\section{Introduction}

In this paper we study the partial regularity of minimizers of the integral functional
\begin{equation}\label{functional}
\F({\bf u}) \colon =\int_{\Omega} f(x,{\bf u},D{\bf u})\,\dd x,
\end{equation}
where $\Omega\subseteq\mathbb{R}^{n}$ is {an open} bounded set and $\bfu : \Omega\rightarrow\mathbb{R}^{N}$, with $n$, $N\ge2$ -- i.e., we consider vectorial minimizers of $\F$.  The growth conditions we impose on $f=f(x,{\bf u},{\bf P})$ are quite general, being as they permit \lq \lq  general growth conditions" with respect to the gradient variable.  This allows us to treat in a unified way  the degenerate (when $p>2$) or singular (when $p<2$) behaviour.  We assume with respect to $x$ a weak VMO condition, uniformly in $(\bfu,\bfP)$, and continuity with respect to $\bfu$.  Our main result, Theorem \ref{theorem-result-1}, proves that a minimizer of \eqref{functional} is locally H\"{o}lder continuous for any H\"{o}lder exponent $0<\alpha<1$ -- i.e., if $\bfu$ is a minimizer of \eqref{functional}, then $\bfu\in C_{\text{loc}}^{0,\alpha}\left(\Omega_0,\mathbb{R}^{N}\right)$, where $\Omega_0\subset\Omega$ is an open set of full measure specified in the statement of Theorem \ref{theorem-result-1} later in this section.

\subsection{Literature Review}

We begin by explaining how the study of functional \eqref{functional} fits into the broader regularity theory research over the past many years.  Before proceeding further, we point out that Mingione \cite{Mingione06Dark} has provided a comprehensive account of the various areas of study within regularity theory for integral functionals and PDEs; it is an excellent reference for those wishing to read a broad overview of the various areas of interest within the larger realm of regularity theory.

As already mentioned we allow $f$ to satisfy a VMO-type condition with respect to $x$.  More precisely the partial map $\displaystyle x\mapsto\frac{f(x,\bfu,\bfP)}{\varphi(|\bfP|)}$ satisfies a uniform VMO condition; here $\varphi$ is an $N$-function -- see condition \ref{ass-4f} later in this section for the precise formulation.  As a consequence we allow a certain controlled discontinuous behavior with respect to the spatial variable in the integrand of \eqref{functional}.
We prove partial H\"older continuity for the local minimizers. The first paper who considered low order regularity (for variational integrals) was the one by Foss \& Mingione \cite{fossmingione1}, where they were assuming continuity with respect to $x$ and $\bfu$. Thereafter Kristensen \& Mingione \cite{KRIMIN06} proved H\"older continuity for convex integral functionals with continuous coefficients for a fixed H\"older exponent depending on the dimension and the growth exponent. Stronger assumptions as Dini-type conditions  \cite{DUGAMIN04}  lead to partial $C^1$-regularity.  It is worth mentioning the uniform porosity of the singular set for Lipschitzian minimizers of quasiconvex functionals, \cite{KRIMIN07}.

The space of  functions with vanishing mean oscillation  (VMO) has been introduced by Sarason in the realm of harmonic analysis, see \cite{Sar75}. It has had several applications in connection with Hardy spaces, Riesz transforms or nonlinear commutators, see \cite{Stein}, \cite{IwMar} and references therein. 
In the early 90's Chiarenza, Frasca and Longo  \cite{chiafralo} studied non-divergence form equations with VMO coefficients by means of singular integrals operators, see also \cite{difazio1}, \cite{difazio2}.

The study of functionals with VMO-type coefficients has been broadened considerably over the past couple decades, see \cite{ragusa2}, \cite{danecek1}.  Recently,  B\"{o}gelein, Duzaar, Habermann, and Scheven \cite{BDHS} considered a functional of the form \eqref{functional} under the assumption that $(x,\bfu,\bfP)\mapsto f(x,\bfu,\bfP)$ satisfies a type of VMO assumption in $x$, uniformly with respect to $\bfu$ and $\bfP$; they further considered an analogous elliptic system of the form $\nabla\cdot a(x,\bfu,D\bfu)=0$, in which, again, the coefficient $a$ was assumed to satisfy a VMO-type condition with respect to its spatial coordinate.  Moreover, the integral functional they studied was assumed to be quasi-convex.  However, unlike our study, they assumed that the growth of $f$ with respect to $\bfP$ was standard $p$-growth, $p\ge 2.$

Similarly, B\"{o}gelein \cite{Bogelein} studied quasi-convex integral functionals in the vectorial case.  But the assumed growth of the integrand with respect to the gradient was standard $p$-growth.  It was also assumed that the map $\displaystyle x\mapsto\frac{f(x,\bfu,\bfP)}{(1+|\bfP|)^p}$ was VMO, uniformly with respect to $\bfu$ and $\bfP$.  B\"{o}gelein, Duzaar, Habermann, and Scheven \cite{bogeleinduzaar2} made some similar assumptions when considering a system of PDEs involving the symmetric part of the gradient $D\bfu$, wherein the coefficients on the symmetric part are VMO.

Goodrich \cite{goodrich1} then further generalized, in part, the results of 
\cite{BDHS} by considering \eqref{functional} in the case where $x\mapsto f(x,\bfu,\bfP)$ was VMO, uniformly with respect to $\bfu$ and $\bfP$, and, furthermore, in which $f$ was only asymptotically convex.

Next, the study of problems with general growth conditions has been initiated by Marcellini in a list of papers \cite{MARCELLINI93,MARCELLINI96, MARPAPI06} and it is now very rich -- see, e.g., \cite{DIESTROVER09,DIESTROVER10,DIELENSTROVER12,CeladaOk,cristiana, Stroffo2020}. In particular, Marcellini \& Papi proved the Lipschitz bound for a solution of an elliptic system with  general growth of Uhlenbeck type. In view of comparison estimates, it is worth mentioning the paper \cite{DIESTROVER09}, where the $C^{1,\alpha}$ regularity is proven via an excess decay estimate. 
Very recently, DeFilippis \& Mingione have relaxed the hypotheses by considering also growth of exponential type (no $\Delta_2$-condition), \cite{cristianarosario}.  

So, we see that many papers in recent years have treated \emph{either} VMO-type coefficient problems \emph{or} general growth problems.  To our knowledge, it seems that the combination of these two generalities has not been considered as we do in this paper.  Thus, the results of this paper significantly generalize many of the previously mentioned papers.

\subsection{Strategy of the proof}
We briefly explain the strategy of the proof of the main result. As a major difficulty with respect to the proof by B\"ogelein or Duzaar et al. in the $p$-setting, we can't rely on homogeneity of the function $\phi$. In particular, an analog of the Campanato excess
\begin{equation*}
\Psi_\alpha(x_0,\varrho):=\varrho^{-\alpha p}\dashint_{B_\varrho(x_0)}|\bfu - (\bfu)_{x_0,\varrho}|^p\,\mathrm{d}x
\end{equation*}
defined there and playing a key role in the iteration process could not be easily handled in the Orlicz setting.

Our strategy is to find carefully the two quantities which play the role both in the non-degenerate and in the degenerate cases. The first leading quantity is the excess functional
\begin{equation*}
\Phi(x_0,\varrho):=\dashint_{B_\varrho(x_0)}\varphi_{|({D{\bf u}})_{x_0,\varrho}|}(|D{\bf u}-({D{\bf u}})_{x_0,\varrho}|)\,\mathrm{d}x
\end{equation*}
(see \eqref{(3.5)}). 
In the non-degenerate case, when
\begin{equation}
\Phi(x_0,\varrho)\leq \varphi(|(D{\bf u})_{x_0,\varrho}|)\,,
\label{(eqND)}
\end{equation}
we linearize the problem, via the $\mathcal{A}$-harmonic approximation \cite{DIELENSTROVER12}.  This procedure, exploiting assumptions \ref{ass-4f}-\ref{ass-5f} and a freezing technique (with respect to the variables $x$ and $\bfu$) based on the Ekeland variational principle, provides a comparison map which is an almost minimizer of the frozen functional and whose gradient is $L^1$-close to that of the original minimizer (see Lemma~\ref{lem:lemma3.7}). Such comparison map is shown to be approximately $\mathcal{A}$-harmonic, and this property is inherited by the minimizer itself via the comparison estimate. This allows to prove an excess-decay estimate, which, in turn, permits the iteration of the rescaled excess $\frac{\Phi(x_0,\varrho)}{\varphi(|(D{\bf u})_{x_0,\varrho}|)}$ and of a ``Morrey-type'' excess 
\begin{equation*}
\Theta(x_0,\varrho):=\varrho\varphi^{-1}\left(\dashint_{B_{\varrho}(x_0)}\varphi(|D{\bf u}|)\,\mathrm{d}x\right)
\end{equation*}
at each scale. Namely,
there exists $\vartheta\in(0,1)$ such that, if the boundedness conditions 
\begin{equation*}
\frac{\Phi(x_0,\varrho)}{\varphi(|(D{\bf u})_{x_0,\varrho}|)}\leq \varepsilon_* \quad \mbox{ and }\quad \Theta(x_0,\varrho)\leq\delta_*
\end{equation*}
hold on some ball $B_\varrho(x_0)$, then 
\begin{equation*}
\frac{\Phi(x_0,\vartheta^m\varrho)}{{\varphi(|(D{\bf u})_{x_0,\vartheta^m\varrho}|)}}\leq \varepsilon_* \quad \mbox{ and }\quad \Theta(x_0,\vartheta^m\varrho)\leq\delta_*
\end{equation*}
hold for every $m=0,1,\dots.$. Therefore, $\Theta(x_0,\varrho)$ is the adequate excess playing the role of $\Psi_\alpha$ in our setting.

In the degenerate case, when
\begin{equation}
\Phi(x_0,\varrho)\geq \kappa \varphi(|(D{\bf u})_{x_0,\varrho}|)
\label{(eqD)}
\end{equation}
for some $\kappa<1$, we perform a different linearization procedure: the assumption \ref{ass-7f} coupled with an analogous freezing argument as before provides, now, the almost $\phi$-harmonicity of the minimizer via the application of the $\phi$-harmonic approximation \cite{DIESTROVER10} to the comparison map. The corresponding excess improvement implies that if the excess is small at radius $\varrho$ it is also small at some smaller radius $\theta\varrho$, for $\theta<1$. The key point in this iteration process is that the boundedness of both the excess $\Phi$ and the Morrey excess $\Theta$ at some scale $\vartheta\theta^{k_0}\varrho$ (``switching radius'') under assumption \eqref{(eqND)} is satisfied exactly when the degenerate bound \eqref{(eqD)} fails and therefore we can proceed the iteration in the non-degenerate regime. Notice that, if on the one hand $|(D{\bf u})_{x_0,\varrho}|$ might blow up in the iteration since we cannot expect $C^1$-regularity, on the other hand the Morrey excess $\Theta(x_0,\theta^k\varrho)$ stays bounded, exactly as it should be for a $C^{0,\alpha}$-regularity result. In addition, if at level $k_0$ the regime is non-degenerate, the behavior stays non-degenerate at any subsequent level $k\geq k_0$, and the iteration can proceed. The smallness of $\Theta$
at any level ensures H\"{o}lder continuity of $\bfu$ in $x_0$ provided the excess functionals $\Phi$ and $\Theta$ are small at some initial radius
$\varrho$ (actually, this holds in a neighborhood of $x_0$, since these smallness conditions are open). Finally, it is then proven that such a smallness condition on the excesses is indeed satisfied on the complement of the set $\Sigma_1\cup\Sigma_2$ of Theorem~\ref{theorem-result-1}.

\subsection{Assumptions and statement of the main result}
We list here the main assumptions on the integral functional that we are going to study throughout the paper.
We assume that $\varphi:[0,\infty)\to[0,\infty)$ is an $N$-function such that
\begin{enumerate}[font={\normalfont},label={($\varphi$\arabic*)}]
\item $\varphi\in C^1([0,\infty))\cap C^2((0,\infty))$;\label{ass-1phi}
\item $0<\mu_1-1\leq \inf_{t>0}\frac{t\varphi''(t)}{\varphi'(t)}\leq \sup_{t>0}\frac{t\varphi''(t)}{\varphi'(t)}\leq\mu_2-1$, for suitable constants $1<\mu_1\leq\mu_2$. \label{ass-2phi}
\end{enumerate}
We may assume, without loss of generality, that $1<\mu_1<2<\mu_2$.

For the precise notation and definitions, as well as the additional assumptions we will require on $\varphi$, we refer to Section~\ref{sec:prelbas}.

We assume the integrand $f:\Omega\times\R^N\times\mathbb{R}^{N\times n}\to\R$, $f=f(x,{\bf u},{\bf P})$ to be Borel-measurable, such that the partial map ${\bf P}\to f(\cdot,\cdot,{\bf P})\in C^1(\mathbb{R}^{N\times n})\cap C^2(\mathbb{R}^{N\times n}\backslash\{{\bf 0}\})$. We will denote by $Df$ and $D^2f$ the corresponding first and second gradients, respectively, for fixed $x$ and $\bf u$. We require $f$ to comply with the following assumptions:

\begin{enumerate}[font={\normalfont},label={(F\arabic*)}]
\item \emph{coercivity:} there exists $\nu>0$ such that 
\begin{equation*}
\nu \phi(|{\bf P}|)\le f(x,{\bf u},{\bf P})-f(x,{\bf u},{\bf 0})
\end{equation*}
uniformly in $x\in \Omega$ and ${\bf u}\in \R^N$, for every ${\bf P}\in\mathbb{R}^{N\times n}$;\label{ass-1f}
\item \emph{$\phi$-growth conditions with respect to the ${\bf P}$ variable:} there exists a constant $L>0$ such that
\begin{equation*}
|D f(x,{\bf u},{\bf P})|\le L \phi'(|{\bf P}|)\,,\quad |D^2 f(x,{\bf u},{\bf P})|\le L \phi''(|{\bf P}|)\,,
\end{equation*}
uniformly in $x\in \Omega$ and ${\bf u}\in \R^N$, for every ${\bf P}\in\mathbb{R}^{N\times n}$ with $|{\bf P}|\neq0$;\label{ass-2f}
\item\emph{$f$ is degenerate quasiconvex;} i.e., 
\begin{equation*}
\begin{split}
&\int_{B}{f (x, {\bf u}, {\bf P}+D\bm\eta(y)) -f(x, {\bf u}, {\bf P})} \, \dd y \geq \, \nu \, \int_{B}\phi'' (\mu+ \kabs{{\bf P}} + \kabs{D\bm\eta(y)} )\, \kabs{D\bm\eta(y)}^2 \, \dd y,\\
\end{split}
\end{equation*}
for every $x \in \Omega$, ${\bf u} \in \R^N$, every ball $B\subset \Omega$, ${\bf P}\in \mathbb{R}^{N\times n}$ and $\bm\eta\in C^{\infty}_0(B, \mathbb{R}^N)$, $\mu\ge 0$;\label{ass-3f}
\item \emph{the function $x\mapsto f ( x, {\bf u}, {\bf P})/\phi(|{\bf P}|)$ satisfies a VMO-condition, uniformly with respect to $({\bf u}, {\bf P})$:}
\begin{equation*}
\kabs{ f (x, {\bf u}, {\bf P})- (f( \cdot, {\bf u}, {\bf P}))_{x_0, r}} \leq  {v}_{x_0}(x, r) \phi( |{\bf P}|)\,,\quad \mbox{ for all $x\in B_r(x_0)$}
\end{equation*}
 where  $x_0\in\Omega$, $r\in(0,1]$ and ${\bf P}\in\R^{N\times n}$ and ${v}_{x_0}:\R^n\times[0,1]\to[0,2L]$ are bounded functions such that
\begin{equation*}
 \lim_{\varrho\to0}{\mathcal{V}}(\varrho)=0 \,, \mbox{ where } {\mathcal{V}}(\varrho):=\sup_{x_0\in\Omega}\sup_{0<r\leq\varrho} \dashint_{B_r(x_0)}{v}_{x_0}(x,r)\,\mathrm{d}x\,,
\end{equation*}
and
\begin{equation*}
(f( \cdot, {\bf u}, {\bf P}))_{x_0, r}:=\frac{1}{|B_r(x_0)|}\int_{B_r(x_0)} f(x, {\bf u}, {\bf P})\,\mathrm{d}x\,;
\end{equation*}\label{ass-4f}
\item \emph{$f$ is uniformly continuous with respect to the ${\bf u}$ variable;} i.e., 
\begin{equation*}
\kabs{ f (x, {\bf u}, {\bf P})- f (x, {\bf u_0}, {\bf P})}\le L\omega ( |{\bf u}-{\bf u_0}|) \phi(|{\bf P}|)\,, 
\end{equation*}
where $\omega:[0,\infty)\to[0,1]$ is a nondecreasing, concave modulus of continuity; i.e., $\lim_{t\downarrow 0}\omega(t)=\omega(0)=0$.\label{ass-5f}
\item \emph{the second derivatives $D^2f$ are H\"older continuous away from $\bf{0}$ with some exponent  $\beta_0 \in (0,1)$} such that uniformly in $(x, {\bf u})$ and for $0<|{\bf P}|\le \frac12 |{\bf Q}|$
\begin{equation*}
\kabs{ D^2 f (x, {\bf u}, {\bf P})- D^2 f (x, {\bf u}, {\bf P}+{\bf Q}) }  \le   c_0 \,\phi'' (\kabs{{\bf Q}})
			\,\kabs{{\bf Q}}^{-\beta_0} \kabs{{\bf P}}^{\beta_0}\,;
\end{equation*}\label{ass-6f}
\item \emph{the function ${\bf P} \to Df(x,{\bf u}, {\bfP})$ behaves asymptotically at $0$	as the $\phi$-Laplacian;} i.e.,
\begin{equation*}
\lim_{t\to 0} \frac{Df(x,{\bf u}, t {\bfP}) }{\phi'(t)}= {\bfP}\,,
\end{equation*}
uniformly in $\{{\bf P}\in\R^{N\times n}:\,\, |{\bf P}|=1\}$ and uniformly for all $x\in\Omega$ and ${\bf u}\in\R^N$.\label{ass-7f}
\end{enumerate}

Our main regularity result  can be stated as follows.  Note that the definition of $\bf V$ appearing in $\Sigma_1$ can be found in \eqref{eq:defV}.

\begin{theorem}
\label{theorem-result-1}
Let $\Omega \subset \R^n$ be an open bounded domain, $\phi$  a convex function satisfying assumptions {\rm\ref{ass-1phi}}--{\rm\ref{ass-G3new}} and consider a minimizer ${\bf u} \in W^{1,\phi}(\Omega,\R^N)$ to the functional {\rm\eqref{functional}} under the assumptions {\rm\ref{ass-1f}}--{\rm\ref{ass-7f}}. Then there exists an open subset $\Omega_0 \subset \Omega$ such that
\begin{equation*}
 \bfu \in C^{0, \alpha}_{\rm{loc}}\left(\Omega_0 ,\R^{N}\right) \qquad \text{and} \qquad \kabs{\Omega \setminus \Omega_0} \, = \, 0
\end{equation*}
for every $\alpha\in (0,1)$. Moreover, $\Omega \setminus \Omega_0\subset \Sigma_1\cup\Sigma_2$ where
\begin{equation*}
\begin{split}
&\Sigma_1:=\left\{x_0\in\Omega:\,\, \mathop{\lim\inf}_{\varrho\searrow 0}\dashint_{B_\varrho(x_0)}|{\bf V}_{|(D\bfu)_{x_0,\varrho}|}(D\bfu-(D\bfu)_{x_0,\varrho})|^2\,\mathrm{d}x>0\right\}\,,\\
&\Sigma_2:=\left\{x_0\in\Omega:\,\, \mathop{\lim\sup}_{\varrho\searrow 0}|(D\bfu)_{x_0,\varrho}|=+\infty\right\}\,.
\end{split}
\end{equation*} 
\end{theorem}

\section{Preliminaries and basic results}\label{sec:prelbas}

\subsection{Some basic facts on $N$--functions} \label{sec:basicNfunctions}

We recall here some elementary definitions and basic results about Orlicz functions. The following definitions and results can be found, e.g., in \cite{Kras, Kufn, Bennett, Adams}. 

A real-valued function $\phi \colon \R^+_0 \to \R^+_0$ is said to be an \emph{$N$-function} if it is convex and satisfies the following conditions: $\phi(0)=0$, $\varphi$ admits the derivative $\phi'$ and this derivative is right continuous, non-decreasing and satisfies $\phi'(0) = 0$, $\phi'(t)>0$ for $t>0$, and $\lim_{t\to \infty} \phi'(t)=\infty$. 

We say that $\phi$ satisfies the \emph{$\Delta_2$-condition} if there exists $c > 0$ such that for all $t \geq 0$ holds $\phi(2t) \leq c\,
\phi(t)$. We denote the smallest possible such constant by $\Delta_2(\phi)$. Since $\phi(t) \leq \phi(2t)$, the $\Delta_2$-condition is equivalent to $\phi(2t) \sim \phi(t)$, where ``$\sim$'' indicates the equivalence between $N$-functions.\par
%
By $L^\phi$ and $W^{1,\phi}$ we denote the classical Orlicz and
Orlicz-Sobolev spaces, i.\,e.\ $f \in L^\phi$ iff $\displaystyle\int
\phi(|{f}|)\,dx < \infty$ and $f \in W^{1,\phi}$ iff $f, D f
\in L^\phi$. The space $W^{1,\phi}_0(\Omega)$ will denote the closure of $C^\infty_0(\Omega)$ in $W^{1,\phi}(\Omega)$.

We define the function $(\phi')^{-1} \colon \R^+_0 \to \R^+_0$ as
\begin{align*}
  (\phi')^{-1}(t) &:= \sup \{ s \in \R^+_0\,:\,
    \phi'(s) \leq t \} .
\end{align*}
If $\phi'$ is strictly increasing, then $(\phi')^{-1}$ is the inverse
function of $\phi'$.  Then $\phi^\ast \colon \R^+_0 \to
\R^+_0$ with
\begin{align*}
  \phi^\ast(t) &:= \int_0^t (\phi')^{-1}(s)\,ds
\end{align*}
is again an $N$-function and $(\phi^\ast)'(t) =
(\phi')^{-1}(t)$ for $t>0$. $\phi^\ast$ is the Young-Fenchel-Yosida conjugate function of
$\phi$.  Note that $\phi^*(t)= \sup_{a \geq 0} (at - \phi(a))$ and
$(\phi^\ast)^\ast = \phi$. When both $\varphi$ and $\varphi^*$ satisfy $\Delta_2$-condition, by elementary convex analysis it is easy to see that for all $\delta>0$ there exists $c_\delta$ (only depending on $\Delta_2(\phi)$ and $\Delta_2(\phi^\ast)$)
such that for all $t, a \geq 0$ it holds that
\begin{equation*}  \label{eq:young}
  at \leq \delta\, \phi(t) + c_\delta\, \phi^\ast(a)\,.
\end{equation*}

\begin{proposition}\label{prop:properties}
Let $\varphi$ be an $N$-function complying with {\rm\ref{ass-1phi}} and {\rm\ref{ass-2phi}}. Then
\begin{description}
\item[(i)] it holds that
\begin{equation}
\phi'(t) \sim t\,\phi''(t)
\label{eq:phi_pp}
\end{equation}
uniformly in $t > 0$. The constants in~\eqref{eq:phi_pp} are called   the {\em characteristics of~$\phi$};
\item[(ii)] it holds that
\begin{equation*}
\mu_1\leq \inf_{t>0}\frac{t\varphi'(t)}{\varphi(t)}\leq \sup_{t>0}\frac{t\varphi'(t)}{\varphi(t)}\leq\mu_2\,;
\label{(2.1celok)}
\end{equation*}
\item[(iii)] the mappings
\begin{equation*}
t\in(0,+\infty)\to \frac{\varphi'(t)}{t^{\mu_1-1}}\,,\,\, \frac{\varphi(t)}{t^{\mu_1}} \mbox{ \,\, and \,\, } t\in(0,+\infty)\to \frac{\varphi'(t)}{t^{\mu_2-1}}\,,\,\, \frac{\varphi(t)}{t^{\mu_2}}
\end{equation*}
are increasing and decreasing, respectively;
\item[(iv)] as for the functions $\varphi$ and $\varphi'$ applied to multiples of given arguments, the following inequalities hold for every $t\geq0$:
\begin{align*}
& a^{\mu_2}\varphi(t) \leq \varphi(at) \leq  a^{\mu_1}\varphi(t) \mbox{ \,\, and \,\, } a^{\mu_2-1}\varphi'(t) \leq \varphi(at) \leq  a^{\mu_1-1}\varphi'(t) \mbox{ \, if \, } 0<a\leq1\,; \\
& a^{\mu_1}\varphi(t) \leq \varphi(at) \leq  a^{\mu_2}\varphi(t) \mbox{ \,\, and \,\,} a^{\mu_1-1}\varphi'(t) \leq \varphi(at) \leq  a^{\mu_2-1}\varphi'(t) \mbox{\, if \,} a\geq1\,.
\end{align*}
\end{description}
\end{proposition}
In particular, from (iv) it follows that both $\varphi$ and $\varphi^*$  satisfy the $\Delta_2$-condition with constants $\Delta_2(\phi)$ and $\Delta_2(\phi^*)$ determined by $\mu_1$ and $\mu_2$. We will denote by $\Delta_2({\phi, \phi^\ast})$ constants depending on $\Delta_2(\phi)$ and $\Delta_2(\phi^*)$. Moreover, for $t>0$ we have
\begin{equation*}
\phi(t) \sim \phi'(t)\,t\,, \qquad \phi(t) \sim \phi''(t)\,t^2\,,\qquad \phi^\ast\big( \phi'(t) \big) \sim \phi^\ast\big( \phi(t)/t \big)\sim \phi(t)\,.
\label{ineq:phiast_phi_p}
\end{equation*}

We recall also that the following inequalities hold for the inverse function $\varphi^{-1}$:
\begin{align}
a^{\frac{1}{\mu_1}}\varphi^{-1}(t)\leq &\varphi^{-1}(at)\leq a^{\frac{1}{\mu_2}}\varphi^{-1}(t)
\label{(2.3a)} 
\end{align}
for every $t\geq0$ with $0<a\leq1$. The same result holds also for $a\geq1$ by exchanging the role of $\mu_1$ and $\mu_2$.

For given $\phi$ we define the associated $N$-function $\psi$ by
\begin{equation*}
  \label{eq:def_psi}
  \psi'(t) := \sqrt{ \phi'(t)\,t\,}.
\end{equation*}

\noindent
Notice that if $\phi$ satisfies
assumption~\eqref{eq:phi_pp}, then also $\phi^*$, $\psi$, and $\psi^*$
satisfy this assumption.

\noindent
Define $\bfV\,:\, \R^{N\times n} \to \R^{N\times n}$ in
the following way:
\begin{equation}
    \label{eq:defV}
    \bfV(\bfQ)=\psi'(|\bfQ|)\frac{\bfQ}{|\bfQ|}\,.
\end{equation}
It is easy to check that
\begin{equation*}
|{\bfV(\bfQ)}|^2 \sim \phi(|{\bfQ}|)\,,
\end{equation*}
uniformly in $\bfQ \in \R^{N\times n}$.

Another important set of tools are the {\rm shifted $N$-functions}
$\{\phi_a \}_{a \ge 0}$ (see \cite{DIEETT08}). We define for $t\geq0$
\begin{equation*}
  \label{eq:phi_shifted}
  \phi_a(t):= \int _0^t \varphi_a'(s)\, \mathrm{d}s\qquad\text{with }\quad
  \phi'_a(t):=\phi'(a+t)\frac {t}{a+t}.
\end{equation*}
We have the following relations:
\begin{align}
&\phi_a(t) \sim \phi'_a(t)\,t\,; \nonumber \\ 
&\phi_a(t) \sim \phi''(a+t)t^2\sim\frac{\varphi(a+t)}{(a+t)^2}t^2\sim \frac{\varphi'(a+t)}{a+t}t^2\,,\label{(2.6b)}\\
& \phi(a+t)\sim [\phi_a(t)+\phi(a)]\,.\label{(2.6c)}
\end{align}
The families $\{\phi_a \}_{a \ge 0}$ and
$\{(\phi_a)^* \}_{a \ge 0}$ satisfy the $\Delta_2$-condition uniformly in $a \ge 0$. 
The connection between $\bfV$ and $\varphi_{a}$ (see {\cite{DIEETT08}})  is the following:
\begin{equation}\label{eq:equivalence}
   |{ \bfV(\bfP) - \bfV(\bfQ)}|^2 \sim  \phi_{|{\bfP}|}(|{\bfP - \bfQ}|)\,,
\end{equation}
uniformly in $\bfP, \bfQ \in \R^{N\times n}$. 
  
The following lemma (see \cite[Corollary~26]{DieKre08}) deals with the \emph{change of shift} for $N$-functions.

\begin{lemma}\label{lem:changeshift}
Let $\varphi$ be an $N$-function with $\Delta_2(\varphi),\Delta_2(\varphi^*)<\infty$. Then for any $\eta>0$ there exists $c_\eta>0$, depending only on $\eta$ and $\Delta_2(\varphi)$, such that for all ${\bf a}, {\bf b}\in\R^d$ and $t\geq0$
\begin{equation}
\varphi_{|{\bf a}|}(t) \leq c_\eta\varphi_{|{\bf b}|}(t) + \eta \varphi_{|{\bf a}|}(|{\bf a}-{\bf b}|)\,.
\label{(5.4diekreu)}
\end{equation}
\end{lemma}

We define the function ${\bf V}_a : \R^{N\times n}\to\R^{N\times n}$ for $a\geq0$ by
\begin{equation}
{\bf V}_a({\bf Q}):= \sqrt{\varphi'_a(|{\bf Q}|)|{\bf Q}|}\frac{{\bf Q}}{|{\bf Q}|}\,,
\label{Vmu}
\end{equation}
where $\varphi_a$ is the shifted $N$-function of $\varphi$. Since $\varphi_0=\varphi$, we retrieve in \eqref{Vmu} the function ${\bf V}$ for $a=0$. With the following lemma, we list some properties of functions ${\bf V}_a$ which will be useful in the sequel.

\begin{lemma}
Let $a\geq0$ and ${\bf V}_a$ be as above. Then for any ${\bf P}, {\bf Q}\in\R^{N\times n}$ a Young-type inequality holds:
\begin{equation}
{\varphi'_a(|{\bf Q}|)}|{\bf P}| \leq c(|{\bf V}_a({\bf Q})|^2+|{\bf V}_a({\bf P})|^2)\,,
\label{young-in}
\end{equation}
where the constant $c$ depends only on $\Delta_2(\varphi)$.
\end{lemma}

Let ${\bf{P}}_0,{\bf{P}}_1\in \mathbb{R}^{N\times n}$, $\theta\in[0,1]$ and define ${\bf{P}}_{\theta}:=(1-\theta){\bf{P}}_0+\theta{\bf{P}}_1$. Then the following result holds (see \cite[Lemma~20]{DIEETT08}).
\begin{lemma}
\label{technisch-mu}
Let $\varphi$ be a $N$-function with $\Delta_2(\varphi, \varphi^*)<\infty.$ Then uniformly for all ${\bf{P}}_0,{\bf{P}}_1\in \mathbb{R}^{N\times n}$ with $|{\bf{P}}_0|+|{\bf{P}}_1|>0$ holds
\begin{equation*}
 \int_0^1 \frac{\varphi'(|{\bf{P}}_{\theta}|)}{|{\bf{P}}_{\theta}|}\, \mathrm{d}\theta \sim \frac{\varphi'(|{\bf{P}}_0|+|{\bf{P}}_1|)}{|{\bf{P}}_0|+|{\bf{P}}_1|}
\end{equation*}
where the constants only depend on $\Delta_2(\varphi, \varphi^*).$
\end{lemma}
In view of the previous considerations, the same proposition holds true for the shifted functions, uniformly in $a\ge 0$.

From assumption \ref{ass-2f} we can easily infer an upper bound for $f(x,{\bf u},{\bf P})-f(x,{\bf u},{\bf Q})$, uniformly in $x\in \Omega$ and ${\bf u}\in \R^N$, for every ${\bf P}, {\bf Q}\in\mathbb{R}^{N\times n}$; namely,
\begin{equation}
\begin{split}
|f(x,{\bf u},{\bf P})-f(x,{\bf u},{\bf Q})|&\leq |{\bf P}-{\bf Q}|\int_0^1 |Df(x,{\bf u},{\bf P}+t({\bf Q}-{\bf P}))|\,\mathrm{d}t \\ 
&\leq L|{\bf P}-{\bf Q}|\int_0^1\varphi'(|{\bf P}+t({\bf Q}-{\bf P})|)\,\mathrm{d}t \\
& \leq cL \varphi(|{\bf P}|+|{\bf Q}|)\,.
\end{split}
\label{(1.3celok)}
\end{equation}

The following estimate is a consequence of \ref{ass-2f} and Lemma~\ref{technisch-mu} (see \cite[eq. (2.14)]{DIELENSTROVER12}):
\begin{equation}
\begin{split}
|Df(x,{\bf u},{\bf P})-Df(x,{\bf u},{\bf Q})|&\leq c(\varphi,L) \varphi''(|{\bf P}|+|{\bf Q}|)|{\bf P}-{\bf Q}|\\
&\leq c(\varphi,L) \varphi'_{|{\bf P}|}(|{\bf P}-{\bf Q}|)\\
& = c(\varphi,L) \frac{\varphi'(|{\bf P}|+|{\bf P}-{\bf Q}|)}{|{\bf P}|+|{\bf P}-{\bf Q}|}|{\bf P}-{\bf Q}|\,,
\label{(2.9)}
\end{split}
\end{equation}
for every ${\bf P}, {\bf Q}\in\R^{N\times n}$. 

The following version of Sobolev-Poincar\'e inequality can be found in \cite[Lemma~7]{DIEETT08}.
\begin{theorem}\label{thm:sob-poincare}
Let $\varphi$ be an $N$-function with $\Delta_2(\varphi,\varphi^*)<+\infty$. Then there exist numbers $\alpha=\alpha(n,\Delta_2(\varphi,\varphi^*))\in(0,1)$ and $K=K(n,N,\Delta_2(\varphi,\varphi^*))>0$ such that the following holds. If $B\subset \R^n$ is any ball with radius $R$ and ${\bf w}\in W^{1,\varphi}(B,\R^N)$, then
\begin{equation}
\dashint_B \varphi\left(\frac{|{\bf w}-({\bf w})_B|}{R}\right)\,\mathrm{d}x\leq K \left(\dashint_B \varphi^\alpha\left({|D{\bf w}|}\right)\,\mathrm{d}x\right)^\frac{1}{\alpha}\,,
\label{eq:sob-poincare}
\end{equation}
where $\displaystyle({\bf w})_B:=\dashint_B {\bf w}(x)\,\mathrm{d}x$. Moreover, if ${\bf w}\in W^{1,\varphi}_0(B,\R^N)$, then
\begin{equation*}
\dashint_B \varphi\left(\frac{|{\bf w}|}{R}\right)\,\mathrm{d}x\leq K \left(\dashint_B \varphi^\alpha\left({|D{\bf w}|}\right)\,\mathrm{d}x\right)^\frac{1}{\alpha}\,,
\label{eq:sob-poincare2}
\end{equation*}
where $K$ and $\alpha$ have the same dependencies as before.
\end{theorem}

\subsection{Some useful lemmas}

The following lemma, useful in order to re-absorb certain terms, is a variant of the classical \cite[Lemma~6.1]{GIUSTI} (see \cite[Lemma~3.1]{DIELENSTROVER12}).

\begin{lemma}\label{lem:iterationlemma}
Let $\psi$ be an $N$-function with $\psi\in\Delta_2$, let $\varrho>0$ and $h\in L^\psi(B_{\varrho}(x_0))$. Let $g:[r,\varrho]\to\R$ be nonnegative and bounded such that for all $r\leq s<t\leq\varrho$
\begin{equation*}
g(s)\leq\theta g(t) + A \int_{B_t(x_0)}\psi\left(\frac{|h(y)|}{t-s}\right)\,\mathrm{d}y+ \frac{B}{(t-s)^\beta}+C\,,
\end{equation*}
where $A,B,C\geq0$, $\beta>0$ and $\theta\in[0,1)$. Then 
\begin{equation*}
g\left(r\right)\leq c(\theta,\Delta_2(\psi),\beta)\left[A\int_{B_{\varrho}(x_0)}\psi\left(\frac{|h(y)|}{\varrho-r}\right)\,\mathrm{d}y + \frac{B}{(\varrho-r)^\beta}+C\right]\,.
\end{equation*}
\end{lemma}

{The following lemma is useful to derive reverse H\"older estimates. It is a variant of the results by Gehring~\cite{Gehring} and Giaquinta-Modica~\cite[Theorem~6.6]{GIUSTI}.
\begin{lemma}\label{lem:gehring}
Let $B_0\subset\R^n$ be a ball, $f\in L^1(B_0)$, and $g\in L^{\sigma_0}(B_0)$ for some $\sigma_0>1$. Assume that for some $\theta\in(0,1)$, $c_1>0$ and all balls $B$ with $2B\subset B_0$
\begin{equation*}
\dashint_B |f|\,\mathrm{d}x\leq c_1 \left(\dashint_{2B}|f|^\theta\,\mathrm{d}x\right)^{1/\theta} + \dashint_{2B}|g|\,\mathrm{d}x\,.
\end{equation*}
Then there exist $\sigma_1>1$ and $c_2>1$ such that $g\in L^{\sigma_1}_{\rm loc}(B)$ and for all $\sigma_2\in[1,\sigma_1]$
\begin{equation*}
\left(\dashint_{B}|f|^{\sigma_2}\,\mathrm{d}x\right)^{1/{\sigma_2}}\leq c_2 \dashint_{2B}|f|\,\mathrm{d}x + c_2 \left(\dashint_{2B}|g|^{\sigma_2}\,\mathrm{d}x\right)^{1/{\sigma_2}}\,.
\end{equation*}
\end{lemma}}

\subsection{$\mathcal{A}$-harmonic and $\varphi$-harmonic functions}

Let $\mathcal{A}$ be a bilinear form on $\R^{N\times n}$. We say that $\mathcal{A}$ is {\em strongly elliptic in the sense of Legendre-Hadamard} if for all $\bm\xi\in \R^N,\bm\zeta\in\R^{n}$ it holds that 
\begin{equation}
  \kappaA \abs{\bm\xi}^2 \abs{\bm\zeta}^2\leq \langle\mathcal{A}(\bm\xi\otimes\bm\zeta)|(\bm\xi\otimes\bm\zeta)\rangle\leq L_{\mathcal{A}} \abs{\bm\xi}^2 \abs{\bm\zeta}^2
\label{(2.20)}
\end{equation}
for some $L_{\mathcal{A}}\geq \kappaA>0$. 
We say that a Sobolev function $\bfw$ on a ball~$B_\varrho(x_0)$ is
\emph{$\mathcal{A}$-harmonic} on $B_\varrho(x_0)$ if it satisfies $-\divergence (\mathcal{A}D \bfw)=0$ in the sense of distributions; i.e.,
\begin{equation*}
\int_{B_\varrho(x_0)} \langle\mathcal{A}D{\bf w}|D\bm\psi\rangle\,\mathrm{d}x=0\,,\quad \mbox{ for all }\bm\psi\in C^\infty_0(B_\varrho(x_0),\R^N)\,.
\end{equation*}
It is well known from the classical theory (see, e.g.~\cite[Chapter 10]{GIUSTI}) that ${\bf w}$ is smooth in the interior of $B_\varrho(x_0)$, and it satisfies the estimate
\begin{equation}
\sup_{B_{\varrho/2}(x_0)}|D {\bf w}|^2+\varrho^2 \sup_{B_{\varrho/2}(x_0)}|D^2 {\bf w}|^2\leq c(n,N,\nu,L) \dashint_{B_\varrho(x_0)}|D{\bf w}|^2\,\mathrm{d}x\,.
\label{(2.21)}
\end{equation}
Let $\varphi$ be an Orlicz function.  We say that a map ${\bf w}\in W^{1,\varphi}(B_\varrho(x_0),\R^N)$ is \emph{$\varphi$-harmonic} on $B_\varrho(x_0)$ (see \cite{DIESTROVER10}) if and only if
\begin{equation*}
\int_{B_\varrho(x_0)} \left\langle\frac{\varphi'(|D{\bf w}|)}{|D{\bf w}|}D{\bf w}\bigg|D\bm\psi\right\rangle\,\mathrm{d}x=0\,,\quad \mbox{ for all }\bm\psi\in C^\infty_0(B_\varrho(x_0),\R^N)\,.
\end{equation*}

More precisely, $D{\bf w}$ and ${\bf V}(D{\bf w})$ are H\"older continuous due to the following decay estimate, see \cite{DIESTROVER09}.
\begin{proposition}\label{generalgrowth} Let $\phi$  be a convex function complying with {\rm\ref{ass-1phi}}, {\rm\ref{ass-2phi}} and
\begin{enumerate}[font={\normalfont},label={($\varphi$3)}]
\item
\begin{align*}
\hspace{0.5cm} & \phi'' \text{ is H\"older continuous off the diagonal:}
 \hspace{1cm}\\
	& \hspace{2cm}  \left|\phi''(s+t)-\phi''(t)\right|\leq c_0\, \phi''
(t)\, \bigg(
    \frac{\left|s\right|}{t} \bigg)^{\beta_0}\,, \quad \beta_0>0\,, 
\hspace{1cm}\\
	& \text{ for all } t>0 \text{ and } s \in \mathbb{R} \text{ with } \left|s\right| < \frac{1}{2} t.
\end{align*}\label{ass-G3new}
\end{enumerate}

\noindent
Then there exist a constant $c\geq 1$ and an exponent $\gamma_0 \in (0,1)$ depending only on $n,N$ and the characteristics of $\phi$, such that the following statement holds true: whenever ${\bf w} \in W^{1,\phi}(B_{R}(x_0),\mathbb{R}^{N})$ is a weak solution of the system
\begin{equation*}
\, {\rm div}\left( \frac{\phi'(\kabs{D \bfu})}{|D \bfu|} \, D \bfu\right) \, = \, 0 \qquad \text{in }B_{R}(x_0)\,,\
\end{equation*}
 then for every $\tau\in(0,1)$ there hold
\begin{equation*}
\begin{split}
\sup_{B_{\tau R/2}(x_0)} \phi(\kabs{D {\bf w}}) \, &\leq \, c \mI{B_{\tau R}(x_0)} \phi(\kabs{D {\bf w}})\,\mathrm{d}x\,,  \\
	\mI{B_{\tau R}(x_0)} \kabs{\bfV(D{\bf w})-(\bfV(D {\bf w}))_{x_0,{\tau R}}}^2 \, \mathrm{d}x \, &\leq \, c \,\tau^{2 \gamma_0} \, 	\mI{B_R(x_0)} \kabs{\bfV(D{\bf w})-(\bfV(D {\bf w}))_{x_0,R}}^2 \, \mathrm{d}x\,.
\end{split}
\end{equation*}
\end{proposition}

This result can be viewed as the Orlicz version of the milestone theorem of Uhlenbeck \cite{UHLENBECK77} for differential forms solving a $p$-harmonic system, see also \cite{LISA13}.

\subsection{Harmonic type approximation results}

We recall here two different harmonic type approximation results. The first one is the \emph{$\mathcal{A}$-harmonic approximation}: given a Sobolev function $\bfu$ on a ball~$B$, we want to find an $\mathcal{A}$-harmonic function~${\bf w}$ which is ``close'' the function $\bfu$. 
It will be the $\mathcal{A}$-harmonic function with the same boundary values as $\bfu$; i.e., a Sobolev function~${\bf w}$
which satisfies
\begin{equation}
  \label{eq:calA1}
  \begin{cases}
    -\divergence (\mathcal{A} D {\bf w})= 0 &\qquad\text{on $B$}
    \\
    {\bf w}= \bfu & \qquad\text{on $\partial B$}
  \end{cases}
\end{equation}
in the sense of distributions.

Setting ${\bf z} := {\bf w} - \bfu$, then~\eqref{eq:calA1} is equivalent to finding
a Sobolev function~${\bf z}$ which satisfies
\begin{equation}
  \label{eq:calA2}
  \begin{cases}
    -\divergence (\mathcal{A} D {\bf z}) = -\divergence(\mathcal{A}
    D \bfu) &\qquad\text{on $B$}
    \\
    {\bf z}= \bfzero &\qquad\text{on $\partial B$}
  \end{cases}
\end{equation}
in the sense of distributions.

The following $\mathcal{A}$-harmonic approximation result in the setting of Orlicz spaces has been proved in \cite[Theorem~14]{DIELENSTROVER12}.
\begin{theorem}
  \label{thm:Aappr_phi}
  Let $B \subset \subset \Omega$ be a ball with radius~$r_B$ and let
  $\widetilde{B} \subset \Omega$ denote either~$B$ or $2B$. Let
  $\mathcal{A}$ be a strongly elliptic (in the sense of
  Legendre-Hadamard) bilinear form on $\R^{N\times n}$.  Let $\psi$ be an N-function with $\psi \in
  \Delta_2(\psi, \psi^*)$ and let $s>1$. Then for every
  $\epsilon>0$, there exists $\delta>0$ only depending on $n$, $N$,
  $\kappa_A$, $\abs{\mathcal{A}}$, $\Delta_2(\psi,\psi^*)$ and
  $s>1$ such that the following holds.  Let $\bfu \in
  W^{1,\psi}(\widetilde{B},\R^N)$ be {\em almost $\mathcal{A}$-harmonic}
  on~$B$ in the sense that

  \begin{align}
    \label{eq:Aappr_ah}
    \biggabs{\dashint_B \langle\mathcal{A}D \bfu | D \bm\eta\rangle\,\mathrm{d}x}
    \leq \delta \dashint_{\widetilde{B}} \abs{D \bfu}\,\mathrm{d}x
    \norm{D \bm\eta}_{L^\infty(B)}
  \end{align}
  for all $\bm\eta \in C^\infty_0(B,\R^N)$. Then the unique solution ${\bf z}
  \in W^{1, \psi}_0(B,\RN)$ of~\eqref{eq:calA2} satisfies
  \begin{equation*}
 \label{eq:Aappr_est}
    \dashint_B \psi\bigg(\frac{\abs{{\bf z}}}{r_B}\bigg)\,\mathrm{d}x +
    \dashint_B \psi(\abs{D {\bf z}})\,\mathrm{d}x \leq \epsilon
    \left(\bigg(\dashint_{\widetilde{B}} \big(\psi(\abs{D
      \bfu})\big)^s \,\mathrm{d}x\bigg)^{\frac 1s}+\dashint_{\widetilde{B}} \psi(\abs{D {\bf u}})\,\mathrm{d}x\right)\,.
  \end{equation*}
\end{theorem}

\begin{remark}\label{rem:thmmodified}
We will exploit the previous approximation result in a slightly modified version. Indeed, following \cite[Lemma~2.7]{CeladaOk}, under the additional assumption
\begin{equation*}
\dashint_{\tilde{B}} \psi(|D{\bf u}|)\,\mathrm{d}x \leq \left(\dashint_{\tilde{B}} [\psi(|D{\bf u}|)]^s\,\mathrm{d}x\right)^{\frac{1}{s}}\leq \psi(\mu)
\end{equation*}
for some exponent $s>1$ and for a constant $\mu>0$, and \eqref{eq:Aappr_ah} replaced by
  \begin{equation*}
        \biggabs{\dashint_B \langle\mathcal{A}D \bfu | D \bm\eta\rangle\,\mathrm{d}x}
    \leq \delta \mu
    \norm{D \bm\eta}_{L^\infty(B)}\,,
  \end{equation*}
it can be seen with minor changes in the proof that the unique solution ${\bf z}
  \in W^{1, \psi}_0(B,\RN)$ of~\eqref{eq:calA2} satisfies
\begin{equation*}
\dashint_B \psi\bigg(\frac{\abs{{\bf z}}}{r_B}\bigg)\,dx +
    \dashint_B \psi(\abs{D {\bf z}})\,dx \leq \epsilon \psi(\mu)\,.
\label{eq:lemma2.7}
\end{equation*}
\end{remark}

Now, moving on to $\varphi$-harmonic functions, the following \emph{$\varphi$-harmonic approximation} lemma (\cite[Lemma~1.1]{DIESTROVER10}) is the extension to general convex functions of the $p$-harmonic approximation lemma \cite{DUMIN04b}, \cite[Lemma~1]{DUMIN04}, and allows to approximate ``almost $\varphi$-harmonic'' functions by $\varphi$-harmonic functions.

\begin{lemma}\label{lem:phiharmapprox}
Let $\varphi$ satisfy assumption \eqref{eq:phi_pp}. For every $\varepsilon>0$ and $\theta\in(0,1)$ there exists $\delta>0$ which only depends on $\varepsilon$, $\theta$, and the characteristics of $\varphi$ such that the following holds. Let $B\subset \R^n$ be a ball and let $\tilde{B}$ denote either $B$ or $2B$. If ${\bf u}\in W^{1,\varphi}(\tilde{B},\R^N)$ is \emph{almost $\varphi$-harmonic} on a ball $B\subset\R^n$ in the sense that
\begin{equation}
\dashint_B \left\langle\frac{\varphi'(|D{\bf u}|)}{|D{\bf u}|}D{\bf u}\biggl|D\bm\eta\right\rangle\,\mathrm{d}x \leq \delta\left(\dashint_{\tilde{B}}\varphi(|D{\bf u}|)\,\mathrm{d}x+\varphi(\|D\bm\eta\|_\infty)\right)
\label{(23Stroffo)}
\end{equation}
for all $\bm\eta\in C^\infty_0(B,\R^N)$, then the unique $\varphi$-harmonic ${\bf w}\in W^{1,\varphi}(B,\R^N)$ with ${\bf w}={\bf u}$ on $\partial B$ satisfies
\begin{equation}
\left(\dashint_B |{\bf V}(D{\bf u})-{\bf V}(D{\bf w})|^{2\theta}\,\mathrm{d}x\right)^{\frac{1}{\theta}}\leq \varepsilon \dashint_{\tilde{B}}\varphi(|D{\bf u}|)\,\mathrm{d}x\,,
\label{(34Stroffo)}
\end{equation}
where ${\bf V}$ is as in \eqref{eq:defV}. 
\end{lemma}

The estimate \eqref{(34Stroffo)} can be improved when $\varphi(|D{\bf u}|)$ satisfies a reverse H\"older inequality as follows (see \cite[Corollary~2.10]{CeladaOk}).

\begin{lemma}\label{corollary2.10}
Let $B\subset\R^n$ be a ball. Let ${\bf u}\in W^{1,\varphi}(2B,\R^N)$ be such that
\begin{equation*}
\left(\dashint_{B}\varphi^{s_1}(|D{\bf u}|)\,\mathrm{d}x\right)^{\frac{1}{s_1}}\leq \tilde{c}_0 \dashint_{2B}\varphi(|D{\bf u}|)\,\mathrm{d}x
\label{celok2.15}
\end{equation*}
for $s_1>1$ and $\tilde{c}_0>0$. Then for every $\varepsilon\in(0,1)$ there exists $\delta_0=\delta_0(n,N,\mu_1,\mu_2,s_1,\tilde{c}_0,\varepsilon)>0$ such that the following holds: if $\bf u$ is almost $\varphi$-harmonic as in \eqref{(23Stroffo)} with $\delta_0$ in place of $\delta$, then the unique $\varphi$-harmonic function ${\bf w}\in W^{1,\varphi}(B, \R^N)$ such that ${\bf w}={\bf u}$ on $\partial B$ satisfies
\begin{equation*}
\dashint_{B} |{\bf V}(D{\bf u})-{\bf V}(D{\bf w})|^{2}\,\mathrm{d}x\leq \varepsilon \dashint_{2B}\varphi(|D{\bf u}|)\,\mathrm{d}x\,.
\label{(2.16celok)}
\end{equation*}
\end{lemma}

\section{Partial regularity for functionals}

\subsection{Caccioppoli inequalities and higher integrability results}

As usual, the first step in proving a regularity theorem for the minimizers of integral functionals is to establish suitable Caccioppoli-type inequalities.

First, we state a ``zero order'' Caccioppoli inequality. The proof is an adaptation to the $\varphi$-setting of \cite[Lemma~3.1]{Bogelein}, we then omit the details (see also \cite[Theorem~2.4]{CeladaOk}).

\begin{lemma}\label{lem:caccioppoli1}
Let ${\bf u}\in W^{1,\varphi}(\Omega,\R^N)$ be a minimizer of the functional \eqref{functional}, under the assumptions {\rm\ref{ass-1f}-\rm\ref{ass-2f}}. Then, for every ${\bf u_0}\in \R^N$ and $x_0\in\Omega$ and all $0<\varrho<{\rm dist}(x_0,\partial\Omega)$ and $r\in[\varrho/2,\varrho)$ there holds
\begin{equation*}
\dashint_{B_r(x_0)} \varphi(|D{\bf u}|)\,\mathrm{d}x\leq c \dashint_{B_\varrho(x_0)}\varphi\left(\frac{|{\bf u} - {\bf u_0}|}{\varrho-r}\right)\,\mathrm{d}x
\label{eq:caccioppoli1}
\end{equation*}
for some constant $c=c(\varphi, L,\nu)>0$.
\end{lemma}

From Lemma~\ref{lem:caccioppoli1} together with the Sobolev-Poincar\'e inequality (Theorem~\ref{thm:sob-poincare}) and Gehring's Lemma (Lemma~\ref{lem:gehring}), one can infer in a standard way the following higher integrability result (see, e.g.,  \cite[Theorem~2.5]{CeladaOk}).

\begin{lemma}\label{lem:higint}
There exist an exponent $s_0=s_0(n,N,\varphi,L,\nu)>1$ and a constant $c$ depending only on $n,N,\varphi,L,\nu$ such that, if ${\bf u}\in W^{1,\varphi}(\Omega;\R^N)$ is a minimizer of the functional \eqref{functional}, complying with {\rm\ref{ass-1f}}-{\rm\ref{ass-2f}}, then the following holds: for every $s\in(1,s_0]$, for any $x_0\in\Omega$, any radius $0<\varrho<{\rm dist}(x_0,\partial\Omega)$ and $r\in[\varrho/2,\varrho)$, one has 
\begin{equation*}
\dashint_{B_r(x_0)} \varphi^{s}(|D{\bf u}|)\,\mathrm{d}x \leq {c}\left(\frac{\varrho}{\varrho-r}\right)^{n(s-1)}\left(\dashint_{B_\varrho(x_0)}\varphi(|D{\bf u}|)\,\mathrm{d}x\right)^{s}\,.
\label{highint0}
\end{equation*}
\end{lemma}

Another useful tool will be the following global higher integrability result on balls for minimizers of \eqref{functional}, {which has been proven in the Orlicz setting for more general integrands in \cite[Lemma~4.3]{cristiana}.}

\begin{lemma}\label{lemma3.3}
Let ${\bf u}\in W^{1,\varphi}(B_r(x_0),\R^N)$ be such that $\varphi(|D{\bf u}|)\in L^{s_0}(B_r(x_0), \R^N)$ for some $s_0>1$. Then there exists an exponent $s=s(n,N,\varphi,L,\nu,s_0)\in(1,s_0]$ and a constant $c=c(n,N,\varphi,L,\nu)$ such that, if ${\bf v}\in {\bf u}+W^{1,\varphi}_0(B_r(x_0),\R^N)$ is a minimizer of the functional $\displaystyle \mathcal{G}[{\bf v}]:=\int_{B_r(x_0)}g(D{\bf v})\,\mathrm{d}x$ with a $C^1$-integrand $g:\R^{N\times n}\to\R$ complying with the growth assumptions
\begin{equation*}
\nu \varphi(|\bm\xi|)\leq g(\bm\xi) \leq L \varphi(1+|\bm\xi|) \quad \mbox{ and }\quad |Dg(\bm\xi)|\leq L \varphi'(|\bm\xi|)
\end{equation*}
for all $\bm\xi\in\R^{nN}$, then we have $\varphi(|D{\bf v}|)\in L^s(B_r(x_0),\R^N)$ and 
\begin{equation*}
\left(\dashint_{B_r(x_0)} \varphi^s(|D{\bf v}|)\,\mathrm{d}x\right)^{\frac{1}{s}}\leq 
 c\left(\dashint_{B_r(x_0)} \varphi^{s_0}(|D{\bf u}|)\,\mathrm{d}x\right)^\frac{1}{s_0}\,.
\end{equation*}
\end{lemma}

We have the following Caccioppoli inequality of second type for local minimizers of \eqref{functional}, involving affine functions.

\begin{lemma}\label{lem:lemma3.4}
There exists a constant $c=c(n,N,\Delta_2(\varphi),\nu,L)>0$ such that, if ${\bf u}\in W^{1,\varphi}(\Omega;\R^N)$ is a minimizer of the functional \eqref{functional} under the assumptions {\rm\ref{ass-1f}}-{\rm\ref{ass-7f}}, and $\bm\ell:\R^n\to\R^N$ is an affine function, say $\bm\ell(x):={\bf u}_0+{\bf Q}(x-x_0)$ for some ${\bf u}_0\in\R^N$ and ${\bf Q}\in\R^{N\times n}$, then for any ball $B_\varrho(x_0)\subseteq\Omega$ with $\varrho\leq\varrho_0$ there holds
\begin{equation*}
\begin{split}
&\dashint_{B_{\varrho/2}(x_0)} \varphi_{|{\bf Q}|}(|D{\bf u}-{\bf Q}|)\,\mathrm{d}x \\
&\leq c \dashint_{B_{\varrho}(x_0)} \varphi_{|{\bf Q}|}\left(\frac{|{\bf u}-{\bm\ell}|}{\varrho}\right)\,\mathrm{d}x + c\varphi(|{\bf Q}|)\left[\omega\left(\dashint_{B_{\varrho}(x_0)}|{\bf u}-{\bf u}_0|+|{\bf u}-\bm\ell|\,\mathrm{d}x\right)^{1-\frac{1}{s}}+ [{\mathcal{V}}(\varrho)]^{1-\frac{1}{s}}\right]
\end{split}
\label{caccioppoliII}
\end{equation*} 
for every $s\in(1,s_0]$ where $s_0$ is that of Lemma~\ref{lem:higint}.
\end{lemma}

\proof
We follow the argument of \cite[Lemma~3.5]{Bogelein} for functionals with $p$-growth, just mentioning how to obtain the analogous main estimates therein. We assume, without loss of generality, that $x_0=0$. For radii $\frac{\varrho}{2}\leq r<\tau<t\leq \frac{3\varrho}{4}$ with $\tau:=\frac{r+t}{2}$ we consider a cut-off function $\eta\in C_0^\infty(B_\tau;[0,1])$ such that $\eta\equiv1$ on $B_r$ and $|D\eta|\leq \frac{4}{t-r}$ on $B_\tau$. Correspondingly, we define the functions $\bm\xi:=\eta({\bf u} - \bm \ell)\in W^{1,\varphi}(B_\tau;\R^N)$ and $\bm\psi:=(1-\eta)({\bf u} - \bm \ell)\in W^{1,\varphi}(B_\tau;\R^N)$. Note that $\bm\ell+\bm\xi={\bf u}-\bm\psi$. From the quasi-convexity assumption \ref{ass-3f}, \eqref{(2.6b)} and simple manipulations we obtain
\begin{equation}
\begin{split}
\int_{B_\tau} \varphi_{|{\bf Q}|}(|D\bm\xi|)\,\mathrm{d}x &\leq c(\nu,\varphi)\int_{B_\tau} \varphi''(|{\bf Q}|+|D\bm\xi|)|D\bm\xi|^2\,\mathrm{d}x \\
&\leq c\int_{B_\tau}[(f(\cdot, {\bf u}_0, {\bf Q}+D\bm\xi(x)))_\tau - (f(\cdot, {\bf u}_0, {\bf Q}))_\tau]\,\mathrm{d}x \\
& = c(J_1+J_2+J_3+J_4+J_5+J_6+J_7)\,,
\end{split}
\label{(3.2)}
\end{equation}
where
\begin{equation*}
\begin{split}
J_1 & := \int_{B_\tau}[(f(\cdot, {\bf u}_0, D{\bf u}(x)-D\bm\psi(x)))_\tau - (f(\cdot, {\bf u}_0, D{\bf u}(x)))_\tau]\,\mathrm{d}x\,,\\
J_2 & := \int_{B_\tau}[(f(\cdot, {\bf u}_0, D{\bf u}(x)))_\tau - (f(\cdot, {\bf u}(x), D{\bf u}(x)))_\tau]\,\mathrm{d}x\,, \\
J_3 & := \int_{B_\tau}[(f(\cdot, {\bf u}(x), D{\bf u}(x)))_\tau - f(x, {\bf u}(x), D{\bf u}(x))]\,\mathrm{d}x\,, \\
J_4 & := \int_{B_\tau}[ f(x, {\bf u}(x), D{\bf u}(x)) - f(x, {\bf u}(x)-\bm\xi(x), D{\bf u}(x)-D\bm\xi(x))]\,\mathrm{d}x\,, \\
J_5 & := \int_{B_\tau}[f(x, {\bf u}(x)-\bm\xi(x), {\bf Q}+D\bm\psi(x))-f(x, {\bf u}_0, {\bf Q}+D\bm\psi(x))]\,\mathrm{d}x\,, \\
J_6 & := \int_{B_\tau}[f(x, {\bf u}_0, {\bf Q}+D\bm\psi(x))-(f(\cdot, {\bf u}_0, {\bf Q}+D\bm\psi(x)))_\tau]\,\mathrm{d}x\,,\\
J_7 & := \int_{B_\tau}[(f(\cdot, {\bf u}_0, {\bf Q}+D\bm\psi(x)))_\tau- (f(\cdot, {\bf u}_0, {\bf Q}))_\tau]\,\mathrm{d}x\,.
\end{split}
\end{equation*}
Now, we proceed to estimate each term above separately. From the minimizing property of ${\bf u}$ we infer that $J_4\leq0$, and by assumptions \ref{ass-5f} and \ref{ass-4f} we obtain the estimates
\begin{equation*}
\begin{split}
J_2 &\leq \int_{B_\tau} \omega(|{\bf u}-{\bf u}_0|)\varphi(|D{\bf u}|)\,\mathrm{d}x\,,\\
J_3 &\leq \int_{B_\tau} {v}_0(\cdot,\tau)\varphi(|D{\bf u}|)\,\mathrm{d}x\,,
\end{split}
\end{equation*}
respectively. Again by exploiting property \ref{ass-5f}, the monotonicity of $\omega$ and $\varphi$, and the fact that 
\begin{equation*}
|\bm\xi|\leq|{\bf u}-\bm\ell|\,, \mbox{\,\, and \,\,} |D\bm\psi|\leq |D{\bf u}-{\bf Q}|+4\left|\frac{{\bf u}-\bm\ell}{t-\tau}\right|\,,
\end{equation*}
we can estimate $J_5$ as
\begin{equation*}
\begin{split}
J_5 & \leq c(\varphi)\int_{B_\tau} \omega(|{\bf u}-{\bf u}_0|+|{\bf u}-\bm\ell|)\varphi\left(|{\bf Q}|+|D{\bf u}|+\left|\frac{{\bf u}-\bm\ell}{t-\tau}\right|\right)\,\mathrm{d}x\\
&\leq c(\varphi) \int_{B_\tau} \omega(|{\bf u}-{\bf u}_0|+|{\bf u}-\bm\ell|)\left[\varphi(|{\bf Q}|+|D{\bf u}|)+\varphi\left(\left|\frac{{\bf u}-\bm\ell}{t-\tau}\right|\right)\right]\,\mathrm{d}x\,,
\end{split}
\end{equation*}
whence, taking into account that by virtue of \eqref{(2.6c)},
\begin{equation}
\begin{split}
\varphi\left(\left|\frac{{\bf u}-\bm\ell}{t-\tau}\right|\right)&\leq c\varphi_{|{\bf Q}|}\left(\left|\frac{{\bf u}-\bm\ell}{t-\tau}\right|\right)+c\varphi(|{\bf Q}|)\\
& \leq c(\varphi)\left|{\bf V}_{|{\bf Q}|}\left(\left|\frac{{\bf u}-\bm\ell}{t-\tau}\right|\right)\right|^2+c\varphi(|{\bf Q}|)
\end{split}
\label{(3.3)}
\end{equation}
and recalling that $\omega\leq1$, we get
\begin{equation*}
J_5\leq c(\varphi) \left(\int_{B_\tau} \left|{\bf V}_{|{\bf Q}|}\left(\left|\frac{{\bf u}-\bm\ell}{t-\tau}\right|\right)\right|^2\,\mathrm{d}x + \int_{B_\tau} \omega(|{\bf u}-{\bf u}_0|+|{\bf u}-\bm\ell|)\varphi\left(|{\bf Q}|+|D{\bf u}|\right)\,\mathrm{d}x\right)\,.
\end{equation*}
For what concerns $J_6$, an analogous computation as for the estimate of $J_5$ based on \eqref{(3.3)} and the VMO assumption \ref{ass-4f} gives
\begin{equation*}
\begin{split}
J_6 & \leq \int_{B_\tau} {v}_0(\cdot,\tau)\varphi(|{\bf Q}+D\bm\psi|)\,\mathrm{d}x \\
 & \leq c(\Delta_2(\varphi)) \left(\int_{B_\tau} \left|{\bf V}_{|{\bf Q}|}\left(\left|\frac{{\bf u}-\bm\ell}{t-\tau}\right|\right)\right|^2\,\mathrm{d}x + \int_{B_\tau} {v}_0(\cdot,\tau)\varphi\left(|{\bf Q}|+|D{\bf u}|\right)\,\mathrm{d}x\right)\,.
\end{split}
\end{equation*}
The terms $J_1$ and $J_7$ can be combined together as
\begin{equation*}
\begin{split}
J_7+J_1=&\int_{B_\tau}\dashint_{B_\tau}\int_0^1\langle Df(y,{\bf u}_0,{\bf Q}+\theta D\bm\psi(x)) - Df(y,{\bf u}_0,{\bf Q})|D\bm\psi(x)\rangle\,\mathrm{d}\theta\mathrm{d}y\mathrm{d}x \\
  & + \int_{B_\tau}\dashint_{B_\tau}\int_0^1\langle Df(y,{\bf u}_0,{\bf Q}) - Df(y,{\bf u}_0,D{\bf u}(x)-(1-\theta)D\bm\psi(x))| D\bm\psi(x)\rangle\,\mathrm{d}\theta\mathrm{d}y\mathrm{d}x\\
&=: J'_7+J'_1\,.
\end{split}
\end{equation*}
From the Cauchy-Schwarz inequality, \eqref{(2.9)} and the fact that $D\bm\psi={\bf 0}$ on $B_r$ we infer
\begin{equation*}
\begin{split}
J'_7 &\leq \int_{B_\tau}\dashint_{B_\tau}\int_0^1|Df(y,{\bf u}_0,{\bf Q}+\theta D\bm\psi(x)) - Df(y,{\bf u}_0,{\bf Q})||D\bm\psi(x)|\,\mathrm{d}\theta\mathrm{d}y\mathrm{d}x \\
&\leq c\int_{B_\tau}\int_0^1\varphi'_{|{\bf Q}|}(\theta |D\bm\psi(x)|)|D\bm\psi(x)|\,\mathrm{d}\theta\mathrm{d}x \\
& \leq c\int_{B_\tau}\varphi_{|{\bf Q}|}(|D\bm\psi(x)|)\,\mathrm{d}x \leq c\int_{B_\tau\backslash B_r}|{\bf V}_{|{\bf Q}|}(D\bm\psi(x))|^2\,\mathrm{d}x\,.
\end{split}
\end{equation*}
We can estimate $J'_1$ analogously, by recalling that $D{\bf u}-(1-\theta)D\bm\psi={\bf Q}+D\bm\xi+\theta D\bm\psi$, $D\bm\psi={\bf 0}$ on $B_r$ and applying the triangle inequality for $\varphi'_{|{\bf Q}|}$, \eqref{(2.9)} and the Young's inequality \eqref{young-in}. 
In this way we get
\begin{equation*}
\begin{split}
J'_1 &\leq \int_{B_\tau}\dashint_{B_\tau}\int_0^1|Df(y,{\bf u}_0,{\bf Q}) - Df(y,{\bf u}_0,D{\bf u}(x)-(1-\theta)D\bm\psi(x))||D\bm\psi(x)|\,\mathrm{d}\theta\mathrm{d}y\mathrm{d}x \\
&\leq c\int_{B_\tau}\int_0^1\varphi'_{|{\bf Q}|}(|D\bm\xi+\theta D\bm\psi|)|D\bm\psi|\,\mathrm{d}\theta\mathrm{d}x \\
& \leq c \int_{B_\tau}\varphi'_{|{\bf Q}|}(|D\bm\psi|)|D\bm\psi|\,\mathrm{d}x + c \int_{B_\tau}\varphi'_{|{\bf Q}|}(|D\bm\xi|)|D\bm\psi|\,\mathrm{d}x\\
& \leq c\int_{B_\tau\backslash B_r}(|{\bf V}_{|{\bf Q}|}(D\bm\psi)|^2+|{\bf V}_{|{\bf Q}|}(D\bm\xi)|^2)\,\mathrm{d}x\,.
\end{split}
\end{equation*}
Recalling the definitions of $\bm\xi$ and $\bm\psi$, by a simple computation we find that
\begin{align*}
D\bm\psi &=(1-\eta)(D{\bf u}-{\bf Q})-D\eta\otimes({\bf u}-\bm\ell)\,, \\
D\bm\xi &=\eta (D{\bf u}-{\bf Q}) + D\eta\otimes({\bf u}-\bm\ell)\,,
\end{align*}
whence
\begin{equation*}
\int_{B_\tau\backslash B_r} (\varphi_{|{\bf Q}|}(|D\bm\psi|) + \varphi_{|{\bf Q}|}(|D\bm\xi|))\,\mathrm{d}x \leq c \int_{B_\tau\backslash B_r} \varphi_{|{\bf Q}|}(|D{\bf u}-{\bf Q}|)\,\mathrm{d}x + c\int_{B_\tau}\varphi_{|{\bf Q}|}\left(\frac{{\bf u}-\bm\ell}{t-r}\right)\,\mathrm{d}x\,, 
\end{equation*}
so that combining with the previous estimates we get
\begin{equation*}
J_1+J_7\leq c\left(\int_{B_\tau\backslash B_r}|{\bf V}_{|{\bf Q}|}(D{\bf u}-{\bf Q})|^2\,\mathrm{d}x+\int_{B_\tau}\left|{\bf V}_{|{\bf Q}|}\left(\frac{{\bf u}-\bm\ell}{t-r}\right)\right|^2\,\mathrm{d}x\right)\,.
\end{equation*}
Since $\bm\xi={\bf u}-\bm\ell$ on $B_r$ and $\tau\leq\varrho$, from \eqref{(3.2)} and the estimates for $J_1-J_7$ we obtain
\begin{equation*}
\begin{split}
&\int_{B_r}|{\bf V}_{|{\bf Q}|}(D{\bf u}-{\bf Q})|^2\,\mathrm{d}x \\
& \leq \tilde{c}\left(\int_{B_\tau\backslash B_r}|{\bf V}_{|{\bf Q}|}(D{\bf u}-{\bf Q})|^2\,\mathrm{d}x+\int_{B_\varrho}\left|{\bf V}_{|{\bf Q}|}\left(\frac{{\bf u}-\bm\ell}{t-r}\right)\right|^2\,\mathrm{d}x\right)\\
& + \tilde{c} \int_{B_\tau} \left(\omega(|{\bf u}-{\bf u}_0|+|{\bf u}-\bm\ell|)+{v}_0(\cdot,\tau)\right)\varphi(|{\bf Q}|+|D{\bf u}|)\,\mathrm{d}x\,.
\end{split}
\end{equation*}
Now, in a standard way we ``fill the hole'' thus obtaining
\begin{equation}
\begin{split}
&\int_{B_r}|{\bf V}_{|{\bf Q}|}(D{\bf u}-{\bf Q})|^2\,\mathrm{d}x \\
& \leq \sigma\int_{B_\tau}|{\bf V}_{|{\bf Q}|}(D{\bf u}-{\bf Q})|^2\,\mathrm{d}x+\int_{B_\varrho}\left|{\bf V}_{|{\bf Q}|}\left(\frac{{\bf u}-\bm\ell}{t-r}\right)\right|^2\,\mathrm{d}x\\
& + \int_{B_\tau} \left(\omega(|{\bf u}-{\bf u}_0|+|{\bf u}-\bm\ell|)+{v}_0(\cdot,\tau)\right)\varphi(|{\bf Q}|+|D{\bf u}|)\,\mathrm{d}x\,,
\end{split}
\label{(3.4)}
\end{equation}
where $\sigma:=\frac{\tilde{c}}{\tilde{c}+1}<1$. In order to bound the latter term further, we exploit the higher integrability result of  Lemma~\ref{lem:higint}. Thus, with fixed $s\in(1,s_0]$, as a consequence of H\"older's inequality, the concavity of $\omega$, the bounds $\omega\leq1$ and ${v}_0\leq 2L$, and Jensen's inequality also we obtain 
\begin{equation*}
\begin{split}
&\int_{B_\tau} \left(\omega(|{\bf u}-{\bf u}_0|+|{\bf u}-\bm\ell|)+{v}_0(\cdot,\tau)\right)\varphi(|{\bf Q}|+|D{\bf u}|)\,\mathrm{d}x\\
& \leq c |B_\tau|\left(\dashint_{B_\varrho}\omega(|{\bf u}-{\bf u}_0|+|{\bf u}-\bm\ell|)^{\frac{s}{s-1}}\,\mathrm{d}x+\dashint_{B_\tau}{v}_0(\cdot,\tau)^{\frac{s}{s-1}}\,\mathrm{d}x\right)^{1-\frac{1}{s}}\left(\dashint_{B_\tau}\varphi^s(|{\bf Q}|)+\varphi^s(|D{\bf u}|)\,\mathrm{d}x\right)^{\frac{1}{s}} \\
& \leq c \tau^n\left(\frac{t}{t-r}\right)^{n(s-1)}\left[\omega\left(\dashint_{B_\varrho}|{\bf u}-{\bf u}_0|+|{\bf u}-\bm\ell|\,\mathrm{d}x\right)^{1-\frac{1}{s}}+{\mathcal{V}}(\tau)^{1-\frac{1}{s}}\right]\dashint_{B_t}\varphi(|{\bf Q}|)+\varphi(|D{\bf u}|)\,\mathrm{d}x \\
& \leq c\left(\frac{\varrho}{t-r}\right)^{n(s-1)}\left[\omega\left(\dashint_{B_\varrho}|{\bf u}-{\bf u}_0|+|{\bf u}-\bm\ell|\,\mathrm{d}x\right)^{1-\frac{1}{s}}+{\mathcal{V}}(\tau)^{1-\frac{1}{s}}\right]\dashint_{B_{3\varrho/4}}\varphi(|{\bf Q}|)+\varphi(|D{\bf u}|)\,\mathrm{d}x\,,
\end{split}
\end{equation*}
where $c=c(n,N,\Delta_2(\varphi),\nu,L)$. This estimate, combined with \eqref{(3.4)} gives
\begin{equation*}
\begin{split}
&\int_{B_r}|{\bf V}_{|{\bf Q}|}(D{\bf u}-{\bf Q})|^2\,\mathrm{d}x \\
& \leq \sigma\int_{B_t}|{\bf V}_{|{\bf Q}|}(D{\bf u}-{\bf Q})|^2\,\mathrm{d}x+c\int_{B_\varrho}\left|{\bf V}_{|{\bf Q}|}\left(\frac{{\bf u}-\bm\ell}{t-r}\right)\right|^2\,\mathrm{d}x\\
& + c\left(\frac{\varrho}{t-r}\right)^{n(s-1)}\left[\omega\left(\dashint_{B_\varrho}|{\bf u}-{\bf u}_0|+|{\bf u}-\bm\ell|\,\mathrm{d}x\right)^{1-\frac{1}{s}}+{\mathcal{V}}(\tau)^{1-\frac{1}{s}}\right]\dashint_{B_{3\varrho/4}}\varphi(|{\bf Q}|)+\varphi(|D{\bf u}|)\,\mathrm{d}x\,, \\
& =: \sigma\int_{B_t}|{\bf V}_{|{\bf Q}|}(D{\bf u}-{\bf Q})|^2\,\mathrm{d}x+c\int_{B_\varrho}\left|{\bf V}_{|{\bf Q}|}\left(\frac{{\bf u}-\bm\ell}{t-r}\right)\right|^2\,\mathrm{d}x + c\left(\frac{\varrho}{t-r}\right)^{n(s-1)}\mathcal{U}\,.
\end{split}
\end{equation*}
Now, since the previous estimate holds for arbitrary radii $r,t$ such that $\varrho/2\leq r<t\leq 3\varrho/4$, the constant $c$ depends only on $n,N,\Delta_2(\varphi),\nu,L$ and $\sigma<1$, as a consequence of Lemma~\ref{lem:iterationlemma} applied with $\beta:=n(s-1)$ we obtain
\begin{equation}
\int_{B_{\varrho/2}}|{\bf V}_{|{\bf Q}|}(D{\bf u}-{\bf Q})|^2\,\mathrm{d}x \leq c \int_{B_\varrho}\left|{\bf V}_{|{\bf Q}|}\left(\frac{{\bf u}-\bm\ell}{\varrho}\right)\right|^2\,\mathrm{d}x + c\mathcal{U}\,.
\label{(3.4bis)}
\end{equation}
In view of Lemma~\ref{lem:caccioppoli1} applied with $\varrho$ in place of $t-s$ and from \eqref{(3.3)} we get
\begin{equation*}
\begin{split}
\int_{B_{3\varrho/4}} \varphi(|D{\bf u}|)\,\mathrm{d}x &\leq c \int_{B_{\varrho}}\varphi\left(\left|\frac{{\bf u}-{\bf u}_0}{\varrho}\right|\right)\,\mathrm{d}x \\
&\leq c  \int_{B_{\varrho}}\varphi\left(\left|\frac{{\bf u}-{\bm\ell}}{\varrho}\right|\right)\,\mathrm{d}x + c \varphi(|{\bf Q}|)\\
& \leq c\left[\int_{B_{\varrho}}\left|{\bf V}_{|{\bf Q}|}\left(\frac{{\bf u}-\bm\ell}{\varrho}\right)\right|^2\,\mathrm{d}x + \varphi(|{\bf Q}|)\right]\,,
\end{split}
\end{equation*}
which combined with \eqref{(3.4bis)} and using the fact that $\omega\leq1$ as well as ${\mathcal{V}}(\varrho)\leq2L$ gives
\begin{equation*}
\begin{split}
\int_{B_{\varrho/2}}|{\bf V}_{|{\bf Q}|}(D{\bf u}-{\bf Q})|^2\,\mathrm{d}x &\leq c \int_{B_\varrho}\left|{\bf V}_{|{\bf Q}|}\left(\frac{{\bf u}-\bm\ell}{\varrho}\right)\right|^2\,\mathrm{d}x  \\
& + c\varrho^n \varphi(|{\bf Q}|)\left[\omega\left(\dashint_{B_\varrho}|{\bf u}-{\bf u}_0|+|{\bf u}-\bm\ell|\,\mathrm{d}x\right)^{1-\frac{1}{s}}+{\mathcal{V}}(\tau)^{1-\frac{1}{s}}\right]\,,
\end{split}
\end{equation*}
where $c=c(n,N,\Delta_2(\varphi),\nu,L)$. The Caccioppoli inequality then follows by taking means on both sides of the latter inequality.
\endproof

We can apply Lemma~\ref{lem:lemma3.4} to affine functions $\bm\ell_{x_0,r}(x):=({\bf u})_{x_0,\varrho}+{\bf Q}(x-x_0)$ for some ${\bf Q}\in\R^{N\times n}$, and the resulting Caccioppoli inequality can be compared with that of \cite[Theorem~3.1]{CeladaOk}. We notice that, apart of an extra VMO term due to assumption \ref{ass-4f}, the dependence of the integrand $f$ also on $\bfu$ implies that the remainder term inside $\omega$; i.e., 
\begin{equation}
R(x_0,\varrho,\bfu,{\bf Q}):=\displaystyle \dashint_{B_{\varrho}(x_0)}|{\bf u}-({\bf u})_{x_0,\varrho}|+|{\bf u}-\bm\ell_{x_0,\varrho}|\,\mathrm{d}x 
\label{eq:remainderca}
\end{equation}
is, in general, non-monotone in the radius $\varrho$. Indeed, it can be estimated from above by the \emph{Morrey-type excess}
\begin{equation}
\Theta(x_0,\varrho):=\varrho\varphi^{-1}\left(\dashint_{B_{\varrho}(x_0)}\varphi(|D{\bf u}|)\,\mathrm{d}x\right)\,,
\label{eq:Theta}
\end{equation}
which fails to be monotone for small $\varrho$ (Lemma~\ref{lem:remainder}(i)). This does not allow, in general, for an application of Gehring's lemma in order to infer an higher integrability result: for this purpose, a suitable ``smallness'' regime \eqref{smallnesscon} has to be imposed (Lemma~\ref{lem:remainder}(ii)).

\begin{lemma}\label{lem:remainder}
Let $\bm\ell_{x_0,\varrho}$ be an affine function as above, and $R(x_0,\varrho,\bfu,\bm\ell_{x_0,\varrho})$ be defined as in \eqref{eq:remainderca}. Then
\begin{enumerate}
\item[\rm (i)] 
\begin{equation}
R(x_0,\varrho,\bfu,{\bf Q}) \leq c \Theta(x_0,\varrho) + \varrho |{\bf Q}|\,.
\label{restostima1}
\end{equation}
In particular, if ${\bf Q}=(D{\bf u})_{x_0,\varrho}$, we have
\begin{equation}
R(x_0,\varrho,\bfu,{\bf Q}) \leq c \Theta(x_0,\varrho) \,.
\label{restostima1bis}
\end{equation}
\item[\rm (ii)] {if the smallness assumption
\begin{equation}
\dashint_{B_{\varrho}(x_0)} \varphi_{|{\bf Q}|}(|D{\bf u}-{\bf Q}|)\,\mathrm{d}x \leq \Lambda\varphi(|{\bf Q}|)
\label{smallnesscon}
\end{equation}
holds for some $\Lambda\in(0,1]$, then there exists a constant $c=c(\varphi)>0$ such that
\begin{equation}
\Bigg( \dashint_{B_{\varrho}(x_0)}|{\bf u}-({\bf u})_{x_0,\varrho}|^{\mu_1}\,\mathrm{d}x \Bigg)^{\frac{1}{\mu_1}}
\leq c\Theta(x_0,\varrho)
\leq c \varrho ( |{\bf Q}|)\,;
\label{restostima2} 
\end{equation}
hence
$$
R(x_0,\varrho,\bfu,{\bf Q}) \leq c \varrho ( |{\bf Q}|)\,.
$$}
\end{enumerate}
\end{lemma}

\proof
(i) First, from Poincar\'e inequality and Jensen's inequality we obtain
\begin{equation*}
\varphi\left(\dashint_{B_\varrho(x_0)} \frac{|{\bf u}-({\bf u})_{x_0,\varrho}|}{\varrho}\,\mathrm{d}x\right) \leq \dashint_{B_\varrho(x_0)} \varphi\left(\frac{|{\bf u}-({\bf u})_{x_0,\varrho}|}{\varrho}\right)\,\mathrm{d}x \leq c\dashint_{B_\varrho(x_0)} \varphi(|D{\bf u}|)\,\mathrm{d}x\,,
\end{equation*}
whence
\begin{equation*}
\dashint_{B_\varrho(x_0)} |{\bf u}-({\bf u})_{x_0,\varrho}|\,\mathrm{d}x \leq c \Theta(x_0,\varrho)\,.
\label{improv1}
\end{equation*}
Then, recalling the definition of $\bm\ell_{x_0,r}$, it is immediate to infer the estimate \eqref{restostima1}. As for \eqref{restostima1bis}, it follows from \eqref{restostima1} since $\varrho|(D{\bf u})_{x_0,\varrho}|\leq c \Theta(x_0,\varrho)$.\\

(ii) We note from \ref{ass-2phi} that $\phi(t^{1/\mu_1})$ is convex for $t\ge 0$. Applying Jensen's inequality, the Poincar\'e type estimate in Theorem~\ref{thm:sob-poincare} and the change-shift formula \eqref{(5.4diekreu)} with ${\bf a}={\bf 0}$, and using  assumption \eqref{smallnesscon}
, we obtain
\begin{equation*}
\begin{split}
\varphi\Bigg(\Bigg(\dashint_{B_\varrho(x_0)} \left[\frac{|{\bf u}-({\bf u})_{x_0,\varrho}|}{\varrho}\right]^{\mu_1}\,\mathrm{d}x\Bigg)^{\frac{1}{\mu_1}}\Bigg)
&\leq \dashint_{B_\varrho(x_0)} \varphi\left(\frac{|{\bf u}-({\bf u})_{x_0,\varrho}|}{\varrho}\right)\,\mathrm{d}x  \leq c\dashint_{B_\varrho(x_0)} \varphi(|D{\bf u}|)\,\mathrm{d}x \\
& \leq c \dashint_{B_{\varrho}(x_0)} \varphi(|D{\bf u}-{\bf Q}|)\,\mathrm{d}x + c \varphi(|{\bf Q}|) \\
& \leq c\dashint_{B_{\varrho}(x_0)} \varphi_{|{\bf Q}|}(|D{\bf u}-{\bf Q}|)\,\mathrm{d}x + c\varphi(|{\bf Q}|)  \leq \varphi(c(|{\bf Q}|))\,,
\end{split}
\end{equation*}
which yields  \eqref{restostima2} up to applying $\varphi^{-1}$ to both sides.
\endproof

Now, we are in position to establish a ``conditioned'' higher integrability result for $\varphi_{|{\bf Q}|}(|D{\bf u}-{\bf Q}|)$, under the smallness assumption \eqref{smallnesscon}. 
The result follows as a consequence of Gehring's lemma with increasing supports (Lemma~\ref{lem:gehring}): 

\begin{corollary}\label{corollary3.2}
If ${\bf u}\in W^{1,\varphi}(\Omega;\R^N)$ is a minimizer of the functional \eqref{functional} under the assumptions {\rm\ref{ass-1f}}-{\rm\ref{ass-7f}}, and ${\bf Q}\in\R^{N\times n}$ is such that \eqref{smallnesscon} holds for some $\Lambda\in(0,1]$, then there exist a  constant $c=c(n,N,\Delta_2(\varphi),\nu, L)>0$ and $\sigma>1$ such that 
\begin{equation}
\begin{split}
&\left(\dashint_{B_{\varrho/2}(x_0)} \varphi_{|{\bf Q}|}^\sigma(|D{\bf u}-{\bf Q}|)\,\mathrm{d}x\right)^{\frac{1}{\sigma}} \\
&\leq c \dashint_{B_{\varrho}(x_0)} \varphi_{|{\bf Q}|}(|D{\bf u}-{\bf Q}|)\,\mathrm{d}x+ c\varphi(|{\bf Q}|)\left[\omega\left(\varrho|{\bf Q}|\right)^{1-\frac{1}{s}}+ [{\mathcal{V}}(\varrho)]^{1-\frac{1}{s}}\right]
\end{split}
\label{eq:caccioppoliIbis}
\end{equation} 
holds for every $s\in(1,s_0]$ where $s_0$ is that of Lemma~\ref{lem:higint}.
\end{corollary}

{\proof

{Let $y\in\Omega$ and $r>0$ be such that $B_{2r}(y)\subset\subset B_\varrho(x_0)$. In view of Lemma~\ref{lem:lemma3.4} applied with $\varrho=2r$, $x_0=y$, $\bfu_0=(\bfu)_{y,2r}$ and an arbitrary ${\bf Q}$, we obtain
\begin{equation}\label{hig2pf1}
\begin{aligned}
&\dashint_{B_{r}(y)} \varphi_{|{\bf Q}|}(|D{\bf u}-{\bf Q}|)\,\mathrm{d}x 
\leq c \dashint_{B_{2r}(y)} \varphi_{|{\bf Q}|}\left(\frac{|\bfu-(\bfu)_{y,2r}-{\bf Q}(x-y)|}{2r}\right)\,\mathrm{d}x\\
&\qquad + c\varphi(|{\bf Q}|)\Bigg[\omega\Bigg(\dashint_{B_{2r}(y)}|\bfu-(\bfu)_{y,2r}|+ |{\bf Q}||x-y|\, \dd x\Bigg)^{1-\frac{1}{s}}+ {\mathcal{V}}(2r)^{1-\frac{1}{s}}\Bigg]\,.
\end{aligned}
\end{equation}
Here, we observe that
\begin{equation*}
\begin{split}
\dashint_{B_{2r}(y)}|\bfu-(\bfu)_{y,2r}|+ |{\bf Q}||x-y|\, \dd x & \le c \dashint_{B_{2r}(y)}|\bfu-(\bfu)_{y,2r}|\, \dd x + c |{\bf Q}|r \\
&\le c \Bigg(\frac{1}{\varrho|{\bf Q}|}\dashint_{B_{2r}(y)}|\bfu-(\bfu)_{y,2r}|\, \dd x+1 \Bigg) \varrho |{\bf Q}|\,,
\end{split}
\end{equation*}
which, recalling that $\omega(ct)\le  c\omega(t)$ when $c\ge 1$ since $\omega$ is concave and $\omega(0)=0$, yields
$$
\omega\Bigg(\dashint_{B_{2r}(y)}|\bfu-(\bfu)_{y,2r}|+|{\bf Q}||x-y|\, \dd x\Bigg)^{1-\frac{1}{s}} 
\le c \Bigg(\frac{1}{\varrho|{\bf Q}|}\dashint_{B_{2r}(y)}|\bfu-(\bfu)_{y,2r}|\, \dd x+1 \Bigg) \omega(\varrho|{\bf Q}|)^{1-\frac{1}{s}}.
$$
Moreover, as $\dashint_{B_{2r}(y)}\bfu-(\bfu)_{y,2r}-{\bf Q}(x-y)\,\dd x={\bf 0}$, by the Sobolev--Poincar\'e type inequality \eqref{eq:sob-poincare},
$$
 \dashint_{B_{2r}(y)} \varphi_{|{\bf Q}|}\left(\frac{|\bfu-(\bfu)_{y,2r}-{\bf Q}(x-y)|}{r}\right)\,\mathrm{d}x \leq c  \Bigg(\dashint_{B_{2r}(y)} \varphi_{|{\bf Q}|}^{\alpha}(|D\bfu-{\bf Q}|)\,\mathrm{d}x\Bigg)^{\frac{1}{\alpha}}  
$$
for some $\alpha\in(0,1)$. Therefore, plugging the preceding two estimates into \eqref{hig2pf1} and taking into account that
\begin{equation*}
\dashint_{B_{2r}(y)}|\bfu-(\bfu)_{y,2r}|\, \dd x \leq 2 \dashint_{B_{2r}(y)}|\bfu-(\bfu)_{x_0,\varrho}|\, \dd x\,,
\end{equation*}
we obtain
$$\begin{aligned}
&\dashint_{B_{r}(y)} \varphi_{|{\bf Q}|}(|D{\bf u}-{\bf Q}|)\,\mathrm{d}x 
\leq c   \Bigg(\dashint_{B_{2r}(y)} \varphi_{|{\bf Q}|}^{\alpha}(|D\bfu-{\bf Q}|)\,\mathrm{d}x\Bigg)^{\frac{1}{\alpha}}\\
&\qquad + c\frac{ \varphi(|{\bf Q}|)\omega(\varrho|{\bf Q}|)^{1-\frac{1}{s}}}{\varrho|{\bf Q}|}\dashint_{B_{2r}(y)}|\bfu-(\bfu)_{x_0,\varrho}|\, \dd x+c\varphi(|{\bf Q}|)\left[\omega\left(\varrho|{\bf Q}|\right)^{1-\frac{1}{s}}+ {\mathcal{V}}(\varrho)^{1-\frac{1}{s}}\right]\,.
\end{aligned}$$
Now, since $\bfu -(\bfu)_{x_0,\varrho}\in L^{\mu_1}(B_{\varrho}(x_0))$, as a consequence of Gehring's lemma there exists $\sigma=\sigma(n,N,\mu_1,\mu_2,\nu,L)\in (1,\mu_1)$ such that
$$
\begin{aligned}
&\Bigg( \dashint_{B_{\rho/2}(x_0)}  \varphi_{|{\bf Q}|}^{\sigma}(|D\bfu -{\bf Q}|)\,\mathrm{d}x \Bigg)^{\frac{1}{\sigma}} 
\leq c   \dashint_{B_{\varrho}(x_0)} \varphi_{|{\bf Q}|}(|D\bfu-{\bf Q}|)\,\mathrm{d}x\\
&\qquad + c\frac{ \varphi(|{\bf Q}|)\omega(\varrho|{\bf Q}|)^{1-\frac{1}{s}}}{\varrho|{\bf Q}|}\Bigg(\dashint_{B_{\varrho}(x_0)}|\bfu-(\bfu)_{x_0,\varrho}|^{\sigma}\, \dd x\Bigg)^{\frac{1}{\sigma}}+c\varphi(|{\bf Q}|)\left[\omega\left(\varrho|{\bf Q}|\right)^{1-\frac{1}{s}}+ {\mathcal{V}}(\varrho)^{1-\frac{1}{s}}\right]\,.
\end{aligned}
$$
Finally, applying Lemma~\ref{lem:remainder} (ii), we obtain \eqref{eq:caccioppoliIbis}.}
\endproof}

We conclude this section by introducing the excess functional and other tools useful in the sequel. Let $\bm L_{x_0,\varrho}:\R^n\to\R^N$ be the affine function associated to ${\bf u}$ defined as
\begin{equation}
\bm L_{x_0,\varrho}(x):=({\bf u})_{x_0,\varrho} + {\bf Q}_{x_0,\varrho}(x-x_0)\,,
\label{eq:faffineL}
\end{equation}
where ${\bf Q}_{x_0,\varrho}:=(D{\bf u})_{x_0,\varrho}$. For $x_0\in\Omega$ and $\varrho\in(0,{\rm dist}(x_0,\partial\Omega))$, $\varrho\leq1$, we define the \emph{excess functional} as
\begin{equation}
\Phi(x_0,\varrho)\equiv\Phi(x_0,\varrho, \bm L_{x_0,\varrho}):=\dashint_{B_\varrho(x_0)}\varphi_{|({D{\bf u}})_{x_0,\varrho}|}(|D{\bf u}-({D{\bf u}})_{x_0,\varrho}|)\,\mathrm{d}x
\label{(3.5)}
\end{equation}
and
\begin{equation}
\Psi(x_0,\varrho):=\dashint_{B_\varrho(x_0)}\varphi\left(\frac{|{\bf u}-({\bf u})_{x_0,\varrho}|}{\varrho}\right)\,\mathrm{d}x\,.
\label{(3.6)}
\end{equation}
Moreover, we define also 
\begin{equation}
H(x_0,\varrho):=\frac{1}{1+(2L)^{1-\frac{1}{s}}}\left([\omega(\varrho|(D\bfu)_{x_0,\varrho}|)]^{1-\frac{1}{s}} + [{\mathcal{V}}(\varrho)]^{1-\frac{1}{s}}\right)\,,
\label{(3.8)}
\end{equation}
and
\begin{equation}
\widetilde{H}(x_0,\varrho):=\frac{1}{1+(2L)^{1-\frac{1}{s}}}\left([\omega(\Theta(x_0,\varrho))]^{1-\frac{1}{s}} + [{\mathcal{V}}(\varrho)]^{1-\frac{1}{s}}\right)\,,
\label{(3.8bis)}
\end{equation}
where $s\in(1,s_0]$ is the exponent of Lemma~\ref{lemma3.3} and $\Theta(x_0,\varrho)$ is the excess defined in \eqref{eq:Theta}. 
Since $\omega\leq 1$ and ${\mathcal{V}}(\varrho)\leq 2L$, we have that $H(x_0,\varrho),\,\widetilde{H}(x_0,\varrho)\leq1$, and 
\begin{equation*}
H(x_0,\varrho)\leq c\widetilde{H}(x_0,\varrho)
\label{comparison1}
\end{equation*}
as a consequence of Lemma~\ref{lem:remainder}(i). Under the smallness assumption $\Phi(x_0,\varrho)\leq \Lambda \varphi(|(D{\bf u})_{x_0,\varrho}|)$, by virtue of Lemma~\ref{lem:remainder}(ii) there exists a constant $\tilde{c}=\tilde{c}(\varphi)$ such that
\begin{equation*}
\frac{1}{\tilde{c}}\widetilde{H}(x_0,\varrho)\leq H(x_0,\varrho)\leq c\widetilde{H}(x_0,\varrho)\,.
\label{comparison2}
\end{equation*}
We can rewrite the Caccioppoli inequality \eqref{eq:caccioppoliIbis} as
\begin{equation}
\Phi(x_0,\varrho/2) \leq c \Phi(x_0,\varrho) + c\varphi(|(D\bfu)_{x_0,\varrho}|)H(x_0,\varrho)\,.
\label{eq:caccioppoliI}
\end{equation} 
Note also that by \eqref{eq:equivalence} and, e.g., \cite[Lemma~A.2]{DiKaSch} we have the following equivalence:
\begin{equation*}
\Phi(x_0,\varrho) \sim \dashint_{B_\varrho(x_0)}|{\bf V}(D{\bf u})-{\bf V}((D{\bf u})_{x_0,\varrho})|^2\,\mathrm{d}x \sim  \dashint_{B_\varrho(x_0)}|{\bf V}(D{\bf u})-({\bf V}(D{\bf u}))_{x_0,\varrho}|^2\,\mathrm{d}x\,.
\label{eq:equivalencebis}
\end{equation*}

In the case $x_0=0$, we will use the shorthands $\Phi(\varrho)$, $\Psi(\varrho)$, $\Theta(\varrho)$, $H(\varrho)$ and $\widetilde{H}(\varrho)$ in place of $\Phi(0,\varrho)$, $\Psi(0,\varrho)$, $\Theta(0,\varrho)$, $H(0,\varrho)$ and $\widetilde{H}(0,\varrho)$, respectively.

\subsection{Comparison maps via Ekeland's variational principle}

The proof of the main results will require suitable comparison functions, which will be constructed with a freezing argument in the variables $(x,{\bf u})$ based on Ekeland's variational principle. We recall below a version of this classical tool, whose proof can be found, e.g., in \cite[Theorem~5.6]{GIUSTI}. 

\begin{lemma}[Ekeland's principle]\label{lem:ekeland}
Let $(X,d)$ be a complete metric space, and assume that $F : X\to[0,\infty]$ be not identically $\infty$ and lower semicontinuous with respect to the metric topology on $X$. If for some $u\in X$ and some $\kappa>0$, there
holds
\begin{equation*}
F(u)\leq \inf_XF + \kappa\,,
\end{equation*}
then there exists $v\in X$ with the properties 
\begin{equation*}
d(u,v)\leq 1 \mbox{\,\, and \,\,} F(v)\leq F(w)+\kappa d(v,w) \quad \forall w\in X\,.
\end{equation*}
\end{lemma}

Although a similar analysis in the Orlicz setting, for integrands $f=f(x,\bm\xi)$, has been performed in \cite[Theorem~3.3]{CeladaOk}, we will follow a quite different argument, which refers to the case of $p$-growth as in \cite[Lemma~3.7]{Bogelein}. We will also specify the appropriate complete metric space $X$, which is not explicitly mentioned in \cite[Theorem~3.3]{CeladaOk}.

To this aim, let $B_\varrho(x_0)\subseteq\Omega$ with $\varrho\leq\varrho_0$ and set
\begin{equation}\label{(3.10)}
g(\bm\xi) \equiv g_{{x}_0,\varrho}(\bm\xi) := (f(\cdot, ({\bf u})_{{x}_0,\varrho},\bm\xi))_{{x}_0,\varrho} \quad \mbox{ for all $\bm\xi\in\R^{N\times n}$,}
\end{equation}
and
\begin{equation}
K(x_0,\varrho):=\widetilde{H}(x_0,\varrho)\Psi(x_0,\varrho)
\label{eq:Krho}
\end{equation}
where $\widetilde{H}(x_0,\varrho)$ and $\Psi(x_0,\varrho)$ are defined as in \eqref{(3.8bis)} and \eqref{(3.6)}, respectively.

 As for the complete metric space $(X,d)$, following \cite[Lemma~4.4]{Ok2} we consider
\begin{equation*}
X:=\left\{{\bf w}\in {\bf u}+W^{1,1}_0(B_{\varrho/2}(x_0)):\,\, \dashint_{B_\varrho/2(x_0)}\varphi(|D{\bf w}|)\,\mathrm{d}x\leq \dashint_{B_\varrho/2(x_0)}\varphi(|D{\bf u}|)\,\mathrm{d}x\right\}
\end{equation*}
with the metric
\begin{equation*}
d({\bf w}_1,{\bf w}_2):= \frac{1}{c_*\varphi^{-1}(K(\varrho))}\dashint_{B_{\varrho/2}(x_0)}|D{\bf w}_1-D{\bf w}_2|\,\mathrm{d}x\,, \quad \mbox{ for } {\bf w}_1,{\bf w}_2\in{\bf u}+W_0^{1,1}(B_{\varrho/2(x_0)},\R^N)\,,
\end{equation*}
and note that the functional 
\begin{equation}
\mathcal{G}[{\bf w}]:= \dashint_{B_{\varrho/2}(x_0)}g(D{\bf w})\,\mathrm{d}x \quad \mbox{ in } {\bf u}+W^{1,1}_0(B_{\varrho/2}(x_0),\R^N)\,,
\label{eq:gfunctional}
\end{equation}
is lower semicontinuous in the metric topology. We would get a comparison map ${\bf v}\in {\bf u} + W^{1,1}_0(B_{\varrho/2}({x}_0),\R^N)$
by proving the following lemma.

\begin{lemma}\label{lem:lemma3.7}
Assume that ${\bf u}\in W^{1,\varphi}(\Omega,\R^N)$ is a minimizer of the functional \eqref{functional}, under the
assumptions {\rm\ref{ass-1f}}-{\rm\ref{ass-6f}}. 
Then there exists a minimizer ${\bf v}\in {\bf u} + W^{1,1}_0(B_{\varrho/2}({x}_0),\R^N)$ of the functional 

\begin{equation*}
\widetilde{\mathcal{G}}[{\bf w}]:= \dashint_{B_{\varrho/2}(x_0)} g(D{\bf w})\,\mathrm{d}x + \frac{K(x_0,\varrho)}{\varphi^{-1}(K(x_0,\varrho))}\dashint_{B_{\varrho/2}(x_0)}|D{\bf v}-D{\bf w}|\,\mathrm{d}x\,,
\end{equation*}
that satisfies 
\begin{equation}
\dashint_{B_{\varrho/2}(x_0)}|D{\bf v}-D{\bf u}|\,\mathrm{d}x \leq c_*\varphi^{-1}(K(x_0,\varrho))
\label{(3.11)}
\end{equation}
for some constant $c_*=c_*(n,N,\Delta_2(\varphi),\nu,L)$. Moreover, ${\bf v}$ fulfills the following Euler-Lagrange variational inequality:
\begin{equation}
\left|\dashint_{B_{\varrho/2}(x_0)}\langle Dg(D{\bf v})|D{\bm\eta}\rangle\,\mathrm{d}x\right|\leq \frac{K(x_0,\varrho)}{\varphi^{-1}(K(x_0,\varrho))}\dashint_{B_{\varrho/2}(x_0)}|D\bm\eta|\,\mathrm{d}x
\label{(3.14)}
\end{equation}
for every $\bm\eta\in C^{\infty}_0(B_{\varrho/2}(x_0),\R^N)$.
\end{lemma}

\proof
We may assume, without loss of generality, that $x_0=0$ and, correspondingly, we use the shorthand $K(\varrho)$ for $K(0,\varrho)$. As a first remark, we recall that from Lemma~\ref{lem:caccioppoli1} with $r=\frac{3}{4}\varrho$ we have
\begin{equation}
\dashint_{B_{3\varrho/4}} \varphi(|D{\bf u}|)\,\mathrm{d}x\leq c \dashint_{B_\varrho}\varphi\left(\frac{|{\bf u} - ({\bf u})_\varrho|}{\varrho}\right)\,\mathrm{d}x = c\Psi(\varrho)\,,
\label{(3.15)}
\end{equation}
where $c=c(\Delta_2(\varphi),L,\nu)$. We then denote by $\tilde{\bf v}\in X$ a minimizer of the functional \eqref{eq:gfunctional}
whose existence is ensured by the direct method under the assumptions \ref{ass-1f}-\ref{ass-2f}. From the minimality of $\tilde{\bf v}$, assumption \ref{ass-1f} and \eqref{(1.3celok)} we get
\begin{equation}
\begin{split}
\dashint_{B_{\varrho/2}}\varphi(|D\tilde{\bf v}|)\,\mathrm{d}x & \leq\frac{1}{\nu}\dashint_{B_{\varrho/2}}g(D\tilde{\bf v})-g({\bf 0})\,\mathrm{d}x\\
& \leq\frac{1}{\nu}\dashint_{B_{\varrho/2}}g(D{\bf u})-g({\bf 0})\,\mathrm{d}x\leq\frac{c(\varphi)L}{\nu}\dashint_{B_{\varrho/2}}\varphi(|D{\bf u}|)\,\mathrm{d}x\,.
\end{split}
\label{eq:(3.16)}
\end{equation}
By the sublinearity of $\varphi$, the Poincar\'e inequality (Theorem~\ref{thm:sob-poincare}), Jensen's inequality and \eqref{(3.15)} this gives
\begin{equation*}
\begin{split}
\varphi\left(\dashint_{B_{\varrho/2}} \frac{|\tilde{\bf v}-({\bf u})_\varrho|}{\varrho}\,\mathrm{d}x\right) & \leq c \left(\dashint_{B_{\varrho/2}} \varphi\left(\frac{|\tilde{\bf v}-{\bf u}|}{\varrho}\right)\,\mathrm{d}x + \dashint_{B_{\varrho/2}} \varphi\left(\frac{|{\bf u}-({\bf u})_\varrho|}{\varrho}\right)\,\mathrm{d}x\right) \\
&  \leq c\left[\left(\dashint_{B_{\varrho/2}} \varphi^\alpha(|D\tilde{\bf v}-D{\bf u}|)\,\mathrm{d}x\right)^{\frac{1}{\alpha}} + \left(\dashint_{B_{\varrho/2}} \varphi^\alpha(|D{\bf u}|)\,\mathrm{d}x\right)^{\frac{1}{\alpha}}\right]  \\
&  \leq c \dashint_{B_{\varrho/2}} \varphi(|D\tilde{\bf v}|) +  \varphi(|D{\bf u}|)\,\mathrm{d}x \\
&\leq c \dashint_{B_{\varrho/2}}\varphi(|D{\bf u}|)\,\mathrm{d}x\,,
\end{split}
\end{equation*}
whence
\begin{equation}
\dashint_{B_{\varrho/2}}|\tilde{\bf v}-({\bf u})_\varrho|\,\mathrm{d}x\leq c \Theta(\varrho)
\label{(3.17)}
\end{equation}
where $c=c(\varphi,n,L,\nu)$. Moreover, as a consequence of the higher integrability results of both Lemma~\ref{lem:higint} and \ref{lemma3.3}, together with \eqref{eq:(3.16)} and \eqref{(3.15)}, we infer the higher integrability result
\begin{equation}
\left(\dashint_{B_{\varrho/2}} \varphi^s(|D\tilde{\bf v}|)\,\mathrm{d}x\right)^{\frac{1}{s}} \leq c \dashint_{B_{3\varrho/4}}\varphi(|D{\bf u}|)\,\mathrm{d}x\leq c\Psi(\varrho)\,,
\label{(3.18)}
\end{equation}
where $c=c(n,N,\varphi,\nu,L)$ and $s=s(n,N,\varphi, \nu, L)\in(1,s_0]$.

Now we prove that ${\bf u}$ is an almost minimizer of the functional $\mathcal{G}$. Indeed, from the minimality of ${\bf u}$ and assumptions \ref{ass-4f}, \ref{ass-3f} we get
\begin{equation*}
\begin{split}
\dashint_{B_{\varrho/2}}f(x,{\bf u},D{\bf u})\,\mathrm{d}x - \mathcal{G}[\tilde{\bf v}] & \leq \dashint_{B_{\varrho/2}}f(x,\tilde{\bf v},D\tilde{\bf v})\,\mathrm{d}x - \mathcal{G}[\tilde{\bf v}] \\
& = \dashint_{B_{\varrho/2}}f(x,\tilde{\bf v},D\tilde{\bf v})-(f(\cdot,\tilde{\bf v},D\tilde{\bf v}))_\varrho\,\mathrm{d}x \\
& + \dashint_{B_{\varrho/2}} ( f(\cdot,\tilde{\bf v},D\tilde{\bf v}))_\varrho - (f(\cdot,({\bf u})_\varrho,D\tilde{\bf v}))_\varrho\,\mathrm{d}x \\
& \leq c(L)\dashint_{B_{\varrho/2}}[{v}_0(\cdot,\varrho)+\omega(|\tilde{\bf v}-({\bf u})_\varrho|)] \varphi(|D\tilde{\bf v}|)\,\mathrm{d}x\,.
\end{split}
\end{equation*}
Then, by using Jensen's inequality, the concavity and sub-linearity of $\omega$, \eqref{(3.17)} and \eqref{(3.18)}, from the previous estimate we obtain
\begin{equation}
\begin{split}
&\dashint_{B_{\varrho/2}}f(x,{\bf u},D{\bf u})\,\mathrm{d}x - \mathcal{G}[\tilde{\bf v}] \\
& \leq c\left[\omega\left(\dashint_{B_{\varrho/2}}|\tilde{\bf v}-({\bf u})_\varrho|\,\mathrm{d}x\right)^{1-\frac{1}{s}}+[{\mathcal{V}}(\varrho)]^{1-\frac{1}{s}}\right]\left(\dashint_{B_{\varrho/2}}\varphi^s(|D\tilde{\bf v}|)\,\mathrm{d}x\right)^{\frac{1}{s}}\\
& \leq c \left[\omega(\Theta(\varrho))^{1-\frac{1}{s}}+[{\mathcal{V}}(\varrho)]^{1-\frac{1}{s}}\right]\Psi(\varrho) = c K(\varrho)\,,
\end{split}
\label{3.18bis}
\end{equation}
where $c=c(n,N,\Delta_2(\varrho), \nu, L)$. Arguing similarly, we can estimate
\begin{equation}
\begin{split}
\mathcal{G}[{\bf u}] & -\dashint_{B_{\varrho/2}}f(x,{\bf u},D{\bf u})\,\mathrm{d}x \\
  & =  \dashint_{B_{\varrho/2}}\left [(f(\cdot,({\bf u})_\varrho,D{\bf u}))_\varrho - f(x,({\bf u})_\varrho,D{\bf u})\right]\,\mathrm{d}x + \dashint_{B_{\varrho/2}} \left[f(x,({\bf u})_\varrho,D{\bf u}) - f(x,{\bf u},D{\bf u})\right]\,\mathrm{d}x \\
& \leq c \left[\omega(\Theta(\varrho))^{1-\frac{1}{s}}+[{\mathcal{V}}(\varrho)]^{1-\frac{1}{s}}\right]\Psi(\varrho) = c K(\varrho)\,,
\end{split}
\label{3.18tris}
\end{equation}
where the constant $c$ has the same dependencies as before. Adding term by term \eqref{3.18bis}-\eqref{3.18tris} and taking into account the minimality of $\tilde{\bf v}$, we infer
\begin{equation*}
\mathcal{G}[{\bf u}]\leq \mathcal{G}[\tilde{\bf v}] + c_*K(\varrho) = \min_{{\bf u}+W_0^{1,1}(B_{\varrho/2},\R^N)} \mathcal{G} + c_*K(\varrho)\,,
\label{hypEkeland}
\end{equation*}
for a constant $c_*=c_*(n, N, \Delta_2(\varphi), \nu, L)$. 
Finally,  Ekeland's variational principle (Lemma~\ref{lem:ekeland}) with the choice $\kappa=c_*K(\varrho)$ provides the existence of  a function ${\bf v}\in X$ with the desired property of minimality for the functional $\widetilde{\mathcal{G}}$ and such that $d({\bf u},{\bf v})\leq1$, which corresponds to \eqref{(3.11)}. The inequality \eqref{(3.14)} follows from the validity of the associated Euler-Lagrange variational inequality for  ${\bf v}$ in a standard way.
\endproof

\subsection{Approximate $\mathcal{A}$-harmonicity and $\varphi$-harmonicity}\label{sec:approxharmon}

In this section, we provide two different linearization strategies for the minimization problem, along the lines of \cite[Section~3.2]{Bogelein}, where an analogous analysis has been performed for functionals with $p$-growth. On the one hand, with Lemma~\ref{lem:lemma3.9} we will show that the minimizer ${\bf u}$ of $\mathcal{F}$ is an almost $\mathcal{A}$-harmonic function for a suitable elliptic bilinear form $\mathcal{A}$. On the other hand, this ${\bf u}$ turns out to be an almost $\varphi$-harmonic function (see Lemma~\ref{lem:lemma4.3}). These results will allow us to apply the $\mathcal{A}$-harmonic approximation lemma, respectively the $\varphi$-harmonic approximation lemma. The proof will require, in both cases, the comparison maps obtained with Lemma~\ref{lem:lemma3.7}.

We start by proving the approximate $\mathcal{A}$-harmonicity of a minimizer to \eqref{functional}. To this aim, only assumptions \ref{ass-1f}-\ref{ass-6f} are required on $f$. 

Let $\bm L_{x_0,\varrho}$ be the affine function associated to ${\bf u}$ as in \eqref{eq:faffineL}, which complies with $\bm L_{x_0,\varrho}(x_0)=({\bf u})_{x_0,\varrho}$ and $D\bm L_{x_0,\varrho}=(D{\bf u})_{x_0,\varrho}=:{\bf Q}_{x_0,\varrho}$.
We set
\begin{equation*}
\mathcal{A}:=\frac{D^2g((D{\bf u})_{x_0,\varrho})}{\varphi''(|(D{\bf u})_{x_0,\varrho}|)}\equiv \frac{\left(D^2f(\cdot, ({\bf u})_{x_0,\varrho},(D{\bf u})_{x_0,\varrho})\right)_{x_0,\varrho}}{\varphi''(|(D{\bf u})_{x_0,\varrho}|)}\,.
\label{(4.1celadaok)}
\end{equation*}

We point out that $\mathcal{A}$ defined above is a bilinear form on $\R^{N\times n}$, satisfying the ellipticity assumption \eqref{(2.20)} by virtue of \ref{ass-2f} and \ref{ass-3f}.

\begin{lemma}\label{lem:lemma3.9}
Let ${\bf u}\in W^{1,\varphi}(\Omega,\R^N)$ be a minimizer of the functional \eqref{functional}, under the assumptions {\rm\ref{ass-1f}}-{\rm\ref{ass-6f}}, and assume that for a ball $B_\varrho(x_0)\subseteq\Omega$ the non-degeneracy assumptions
\begin{equation*}
\Phi(x_0,\varrho)\leq \varphi(|(D{\bf u})_{x_0,\varrho}|)\quad \mbox{ and }\quad \varrho\leq 1\,,
\label{(3.19)}
\end{equation*}
are satisfied. Then, ${\bf u}$ is \emph{approximately $\mathcal{A}$-harmonic} on the ball $B_{\varrho/2}(x_0)$, in the sense that there exists $\beta_1=\beta_1(n,N,\mu_1,\mu_2,\nu,L,\beta_0)\in(0,\frac{1}{2})$ such that
\begin{equation}
\begin{split}
\biggl|\dashint_{B_{\varrho/2}(x_0)} &\langle\mathcal{A}(D{\bf u}-(D{\bf u})_{x_0,\varrho})|D\bm\eta\rangle\,\mathrm{d}x\biggr| \\
& \leq c |(D{\bf u})_{x_0,\varrho}|\|D\bm\eta\|_\infty\left\{[{H}(x_0,\varrho)]^{\beta_1}+\frac{\Phi(x_0,\varrho)}{\varphi(|(D{\bf u})_{x_0,\varrho}|)}+\left(\frac{\Phi(x_0,\varrho)}{\varphi(|(D{\bf u})_{x_0,\varrho}|)}\right)^{\frac{1+\beta_0}{2}}\right\}
\end{split}
\label{stimone}
\end{equation}
holds for every $\bm\eta\in C^\infty_c(B_{\varrho/2}(x_0),\R^N)$ for some constant $c = c(n,N,\mu_1, \mu_2,\nu,c_0, L)>0$, where
$\mu_1,\mu_2$ are the characteristics of $\varphi$ and $c_0,\beta_0$ are the constants of assumption {\rm\ref{ass-6f}}.
\end{lemma}

\proof
See \cite[Lemma~4.1]{CeladaOk}. 
\endproof

If, in addition, $f$ complies also with \ref{ass-7f}, we can show that each local minimizer of the functional $\mathcal{F}({\bf u})$ (eq. \eqref{functional}) is almost $\varphi$-harmonic. 

For this, we preliminarly note (see \cite[eq. (4.19)-(4.20)]{CeladaOk}) that assumption \ref{ass-7f} implies the following: 
\begin{equation}
\mbox{for every $\delta>0$, there exists $\sigma=\sigma(\delta)>0$ such that\,\,}\left|Dg({\bf P})-\frac{{\bf P}}{|{\bf P}|}\varphi'(|{\bf P}|)\right|\leq \delta \varphi'(|{\bf P}|)\,,
\label{(4.20)}
\end{equation}
for every ${\bf P}\in\R^{N\times n}$ with $0<|{\bf P}|\leq \sigma$, where the function $g$ has been introduced in \eqref{(3.10)}.

We then have the following result.

\begin{lemma}\label{lem:lemma4.3}
Let ${\bf u}\in W^{1,\varphi}_{\rm loc}(\Omega,\R^N)$ be a local minimizer of the functional \eqref{functional}, and assume that $f$ complies also with {\rm\ref{ass-7f}}. Then there exists $\beta_2=\beta_2(n,N,\mu_1,\mu_2,c_0,L)\in(0,\frac{1}{2})$ such that, for every $\delta>0$ and for $\sigma=\sigma(\delta)>0$ given by \eqref{(4.20)}, the inequality
\begin{equation*}
\begin{split}
&\left|\dashint_{B_{\varrho/2}(x_0)}\left\langle\frac{\varphi'(|D{\bf u}|)}{|D{\bf u}|}D{\bf u}\bigg|D\bm\eta\right\rangle\,\mathrm{d}x\right| \\
& \leq c\left(\delta+[\widetilde{H}(x_0,\varrho)]^{\beta_2}+\frac{\varphi^{-1}(\Psi(x_0,\varrho))}{\sigma}\right) \left(\dashint_{B_{\varrho}(x_0)}\varphi(|D{\bf u}|)\,\mathrm{d}x+\varphi(\|D\bm\eta\|_\infty)\right)
\end{split}
\end{equation*}
holds for every $\bm\eta\in C^\infty_c(B_{\varrho/2}(x_0),\R^N)$ for some constant $c=c(n,N,\mu_1,\mu_2,c_0,\nu,L)>0$. 
\end{lemma}

\proof
See \cite[Lemma~4.3]{CeladaOk}.
\endproof

\subsection{Excess decay estimates: the non-degenerate regime}\label{sec:nondegenerate}

We start by establishing excess improvement estimates in the \emph{non-degenerate regime} characterized by \eqref{(4.5a)} below, i.e. the fact that $\Phi(x_0,\varrho)\leq c\varphi(|(D{\bf u})_{x_0,\varrho}|)$. The strategy of the proof is to exploit Lemma~\ref{lem:lemma3.9} to approximate the given minimizer by $\mathcal{A}$-harmonic functions, for which suitable decay estimates are available from Theorem~\ref{thm:Aappr_phi}. 

We introduce the \emph{hybrid excess functional}
\begin{equation}
\Phi_*(x_0,\varrho):=\Phi(x_0,\varrho)+\varphi(|({D{\bf u}})_{x_0,\varrho}|)[H(x_0,\varrho)]^{\beta_1}\,,
\label{(4.11)}
\end{equation}
where $\beta_1$ is the exponent of Lemma~\ref{lem:lemma3.9}. 
Since $\beta_1<1/2$ and $H(x_0,\varrho)\leq 1$, we deduce, in particular, that
$H(x_0,\varrho)\leq [H(x_0,\varrho)]^{\beta_1}$. Thus, the Caccioppoli inequality \eqref{eq:caccioppoliI} can be re-read as
\begin{equation*}
\Phi(x_0,\varrho/2)\leq c\Phi_*(x_0,\varrho)\,,
\label{excess0}
\end{equation*}
where $c=c(n,N,\mu_1,\mu_2,\nu,L)$.

\begin{lemma}\label{lem:lemma3.12}
For every $\varepsilon\in(0,1)$ there exist $\delta_1,\delta_2\in(0,1]$, where $\delta_i=\delta_i(n,N,\mu_1,\mu_2, \beta_0, \nu,L,\varepsilon)$, $i=1,2$, with the following property: if
\begin{align}
\frac{\Phi(x_0,\varrho)}{\varphi(|(D{\bf u})_{x_0,\varrho}|)}\leq \delta_1 \label{(4.5a)}\\
[H(x_0,\varrho)]^{\beta_1}\leq \delta_2 \label{(4.5b)}
\end{align}
then the excess improvement estimate
\begin{equation}
\Phi(x_0,\vartheta\varrho)\leq c_{\rm dec}\vartheta^2\left[1+\frac{\varepsilon}{\vartheta^{n+2}}\right]\Phi_*(x_0,\varrho)
\label{(4.6)}
\end{equation}
holds for every $\vartheta\in(0,1)$ for some constant $c_{\rm dec}=c_{\rm dec}(n,N,\mu_1,\mu_2,\nu,L,c_1)>0$, where $\Phi_*$ is defined in \eqref{(4.11)}.
\end{lemma}

\proof
The proof follows the argument of \cite[Lemma~4.2]{CeladaOk}. We emphasize that Corollary~\ref{corollary3.2} is crucial in order to obtain the estimate
\begin{equation*}
\begin{split}
\left(\dashint_{B_{\varrho/2}}\left[\frac{\varphi_{|{\bf Q}_{\varrho}|}(|D{\bf u}-{\bf Q}_{\varrho}|)}{\varphi(|{\bf Q}_{\varrho}|)}\right]^{s_0}\,\mathrm{d}x\right)^{\frac{1}{s_0}} & \leq \frac{c}{\varphi(|{\bf Q}_{\varrho}|)}\dashint_{B_{\varrho}}\varphi_{|{\bf Q}_{\varrho}|}(|D{\bf u}-{\bf Q}_{\varrho}|)\,\mathrm{d}x + c[H(\varrho)]^{\beta_1}\\
& \leq {c} \frac{\Phi_*(\varrho)}{\varphi(|{\bf Q}_{\varrho}|)}\,,
\end{split}
\end{equation*}
which comes into play in applying the $\mathcal{A}$-harmonic approximation theorem in the modified version of  Remark~\ref{rem:thmmodified}.
\endproof

\begin{lemma}
Let $\vartheta\in(0,1)$, and assume that 
\begin{equation}
\frac{\Phi(x_0,\varrho)}{\varphi(|(D{\bf u})_{x_0,\varrho}|)}\leq \frac{\vartheta^n}{2^{\mu_2+1} c_{\mu_2}}\,,
\label{eq:smallnessbis}
\end{equation}
where $c_{\mu_2}$ is the constant of the change of shift formula \eqref{(5.4diekreu)} with $\eta=\frac{1}{2^{\mu_2+1}}$.
Then it holds that
\begin{equation}
|(D\bfu)_{x_0,\varrho}|\leq 2|(D\bfu)_{x_0,\vartheta\varrho}|\,.
\label{eq:meanscomp}
\end{equation}
\label{lem:meanscomp}
\end{lemma}

\proof
As a consequence of  \eqref{(5.4diekreu)} for $\eta=\frac{1}{2^{\mu_2+1}}$ and with \eqref{eq:smallnessbis} we get
\begin{equation*}
\begin{split}
\varphi(|(D\bfu)_{x_0,\varrho}-(D\bfu)_{x_0,\vartheta\varrho}|) & \leq \dashint_{B_{\vartheta\varrho}(x_0)} \varphi(|D\bfu - (D\bfu)_{x_0,\varrho}|)\,\mathrm{d}x \\
& \leq c_{\mu_2}\vartheta^{-n}\Phi(x_0,\varrho) + \frac{1}{2^{\mu_2+1}}\varphi(|(D\bfu)_{x_0,\varrho}|) \\
& \leq \frac{1}{2^{\mu_2}} \varphi(|(D\bfu)_{x_0,\varrho}|)\,,
\end{split}
\label{eq:computation}
\end{equation*}
whence, passing to $\varphi^{-1}$ and taking into account \eqref{(2.3a)}, we obtain
\begin{equation*}
|(D\bfu)_{x_0,\varrho}-(D\bfu)_{x_0,\vartheta\varrho}| \leq \frac{1}{2}|(D\bfu)_{x_0,\varrho}|\,.
\end{equation*}
Now,
\begin{equation*}
\begin{split}
|(D\bfu)_{x_0,\varrho}| \leq |(D\bfu)_{x_0,\varrho}-(D\bfu)_{x_0,\vartheta\varrho}|+|(D\bfu)_{x_0,\vartheta\varrho}| \leq \frac{1}{2}|(D\bfu)_{x_0,\varrho}| + |(D\bfu)_{x_0,\vartheta\varrho}|\,,
\end{split}
\end{equation*}
whence \eqref{eq:meanscomp} follows by re-absorbing the first term of the right-hand side into the left. 
\endproof

The excess-decay estimate \eqref{(4.6)} can be iterated, as the non-degeneracy conditions \eqref{(4.5a)}-\eqref{(4.5b)} are also satisfied on any smaller ball $B_{\vartheta^m\varrho}(x_0)$, $m\in\mathbb{N}$, $\vartheta<1$.

\begin{lemma}
Let $\Phi(x_0,\varrho)$ and  $\Theta(x_0,\varrho)$ 
be defined as in \eqref{(3.5)} and \eqref{eq:Theta}, 
respectively. 
Then there exist constants $\delta_*$, $\varepsilon_*$, $\varrho_*\in(0,1]$ and $\vartheta$ such that the following holds:
if the conditions 
\begin{equation}
\frac{\Phi(x_0,\varrho)}{\varphi(|(D{\bf u})_{x_0,\varrho}|)}\leq \varepsilon_* \quad \mbox{ and }\quad \Theta(x_0,\varrho)\leq\delta_*\,. 
\label{eq:0step}
\end{equation}
hold on $B_\varrho(x_0)\subseteq\Omega$ for $\varrho\in(0,\varrho_*]$, then 
\begin{equation}
\frac{\Phi(x_0,\vartheta^m\varrho)}{{\varphi(|(D{\bf u})_{x_0,\vartheta^m\varrho}|)}}\leq \varepsilon_* \quad \mbox{ and }\quad \Theta(x_0,\vartheta^m\varrho)\leq\delta_*
\label{eq:kstep}
\end{equation}
for every $m=0,1,\dots.$. As a consequence, for any $\alpha\in(0,1)$ the following Morrey-type estimate holds:
\begin{equation}
\Theta(y,r)\leq c\delta_*\left(\frac{r}{\varrho}\right)^\alpha
\label{(5.10Stroffo)}
\end{equation}
for all $y\in B_{\varrho/2}(x_0)$ and $r\in(0,\varrho/2]$.
\label{lem:lemma3.13}
\end{lemma}

\proof
As usual, we omit the explicit dependence on $x_0$. Let $\vartheta\in(0,1)$ be such that
\begin{equation}
\vartheta\leq \min\left\{(6c_{\rm dec}2^{\mu_2})^{-\frac{1}{2}},\frac{1}{2^{\mu_2}},\frac{1}{2^{\frac{\mu_2}{\mu_1(1-\alpha)}}}\right\}\,,
\label{eq:choosetheta}
\end{equation}
where $c_{\rm dec}$ is the constant of Lemma~\ref{lem:lemma3.12} depending only on $n,N,\mu_1,\mu_2,\nu,L,c_0$.
Correspondingly, let $\delta_i=\delta_i(n,N,\mu_1,\mu_2,\beta_0,\nu,L,\vartheta)$, $i=1,2$ be the constants of Lemma~\ref{lem:lemma3.12}, applied with the choice $\varepsilon=\vartheta^{n+2}$. We choose $\varepsilon_*>0$ such that
\begin{equation}
\varepsilon_*\leq\min\left\{\frac{\delta_1}{3}, \frac{\delta_2}{2}, \frac{\vartheta^n}{\max\{2c_{\frac{1}{2}},2^{\mu_2+1} c_{\mu_2}\}}\right\}\,,
\label{(3.46verena)}
\end{equation}
where $c_{\frac{1}{2}}$ is the constant in the change-shift formula \eqref{(5.4diekreu)} with $\eta=\frac{1}{2}$, and we fix the constant $\delta_*>0$ so small that
\begin{equation}
\left(\frac{\omega(\delta_*)^{1-\frac{1}{s}}}{1+(2L)^{1-\frac{1}{s}}}\right)^{\beta_1}<\varepsilon_*\,.
\label{(3.47verena)}
\end{equation}
Moreover, we choose a radius $\varrho_*>0$ such that
\begin{equation}
\varrho_*\leq1 \qquad \mbox{and} \qquad \left(\frac{\mathcal{V}(\varrho_*)^{1-\frac{1}{s}}}{1+(2L)^{1-\frac{1}{s}}}\right)^{\beta_1}<\varepsilon_*\,.
\label{(3.48verena)}
\end{equation}
As a consequence, $\varepsilon_*$, $\delta_*$ and $\varrho_*$ have the same dependencies as $\delta_1$, $\delta_2$. In addition, $\delta_*$ depends also on $\omega$, while $\varrho_*$ also on $\omega$ and $\mathcal{V}$.

We argue by induction on $m$. Since \eqref{eq:kstep} are trivially true for $m=0$ by assumption \eqref{eq:0step}, our aim is to show that if \eqref{eq:kstep} holds for some $m\geq1$, then the corresponding inequalities hold with $m+1$ in place of $m$.
Setting
\begin{equation*}
E(B_{\vartheta^m\varrho}):=\dashint_{B_{\vartheta^{m}\varrho}}\varphi(|D{\bf u}|)\,\mathrm{d}x\,,
\end{equation*}
in order to prove the second inequalities in \eqref{eq:kstep} it will suffice to show that
\begin{equation}
E(B_{\vartheta^m\varrho})\leq \varphi\left(\frac{\delta_*}{\vartheta^m\varrho}\right)\,.
\label{eq:equivalent}
\end{equation}
We have, with \eqref{eq:kstep} at step $m$, the shift-change formula \eqref{(5.4diekreu)} with $\eta=\frac{1}{2}$ and \eqref{(3.46verena)}, the estimate
\begin{equation}
\begin{split}
E(B_{\vartheta^{m+1}\varrho}) & \leq 2^{\mu_2-1}\left(c_{\frac{1}{2}}\vartheta^{-n} \Phi(\vartheta^m\varrho) + \frac{1}{2}\varphi(|(D{\bf u})_{\vartheta^m\varrho}|)+ \varphi(|(D{\bf u})_{\vartheta^m\varrho}|)\right)\\
& \leq 2^{\mu_2-1}\left(c_{\frac{1}{2}}\vartheta^{-n} \Phi(\vartheta^m\varrho) +\frac{3}{2} E(B_{\vartheta^m\varrho})\right) \\
& \leq 2^{\mu_2-1}\left(c_{\frac{1}{2}}\vartheta^{-n}\varepsilon_*+\frac{3}{2}\right) E(B_{\vartheta^m\varrho})\\
&\leq 2^{\mu_2-1}\left(c_{\frac{1}{2}}\vartheta^{-n}\varepsilon_*+\frac{3}{2}\right)\vartheta\varphi\left(\frac{\delta_*}{\vartheta^{m+1}\varrho}\right)\\
& \leq \varphi\left(\frac{\delta_*}{\vartheta^{m+1}\varrho}\right)\,.
\end{split}
\label{estim1}
\end{equation}

Now, we prove by induction the first inequality in \eqref{eq:kstep} for $m+1$.
From \eqref{eq:kstep} at step $k$ and the choices of $\delta_*$ and $\varrho_*$ as in \eqref{(3.47verena)}-\eqref{(3.48verena)}, we have
\begin{align*}
&\frac{\Phi(\vartheta^m\varrho)}{{\varphi(|(D{\bf u})_{\vartheta^m\varrho}|)}}\leq \varepsilon_*<3\varepsilon_*\leq\delta_1\,,\\
&[H(\vartheta^m\varrho)]^{\beta_1}<2\varepsilon_*\leq \delta_2\,,
\end{align*}
and
\begin{equation*}
\begin{split}
\Phi_*(\vartheta^m\varrho)=\Phi(\vartheta^m\varrho) + \varphi(|(D{\bf u})_{\vartheta^m\varrho}|) [H(\vartheta^m\varrho)]^{\beta_1} & \leq 3\epsilon_*\varphi(|(D{\bf u})_{\vartheta^m\varrho}|)\,.
\end{split}
\end{equation*}
Then, by virtue of Lemma~\ref{lem:lemma3.12} and Lemma~\ref{lem:meanscomp} applied with radius $\vartheta^m\varrho$ in place of $\varrho$, and recalling the choice of $\vartheta$ \eqref{eq:choosetheta}, we get
\begin{equation*}
\begin{split}
\Phi(\vartheta^{m+1}\varrho)\leq 2c_{\rm dec}\vartheta^2\Phi_*(\vartheta^m\varrho)&\leq 6c_{\rm dec}\epsilon_*\vartheta^2\varphi(|(D{\bf u})_{\vartheta^m\varrho}|) \\
& \leq \epsilon_* \varphi(|(D{\bf u})_{\vartheta^{m+1}\varrho}|)\,.
\end{split}
\end{equation*}

Finally, since the iteration starting from $m=0$ of the estimate $\varphi^{-1}(E(B_{\vartheta^{m+1}\varrho}))\leq 2^{\frac{\mu_2}{\mu_1}}\varphi^{-1}(E(B_{\vartheta^{m}\varrho}))$, obtained by \eqref{estim1} and \eqref{(2.3a)}, with \eqref{eq:choosetheta} yields 
\begin{equation*}
(\vartheta^m\varrho)^{1-\alpha}\varphi^{-1}(E(B_{\vartheta^m\varrho}))\leq \varrho^{1-\alpha}\varphi^{-1}(E(B_{\varrho}))\leq \delta_*\varrho^{-\alpha}\,,
\end{equation*}
and this estimate \emph{a fortiori} holds if we consider $E(B_{\vartheta^m\varrho}(y))$ for $y\in B_{\varrho/2}$ in place of $E(B_{\vartheta^m\varrho})$, we deduce the Morrey-type estimate
\begin{equation*}
r^{1-\alpha}\varphi^{-1}(E(B_{r}(y)))\leq c\delta_*\varrho^{-\alpha}
\end{equation*}
for all $y\in B_{\varrho/2}$ and $r\leq\varrho/2$, which is equivalent to \eqref{(5.10Stroffo)}. The proof is now concluded.
\endproof

\subsection{Excess decay estimate: the degenerate regime}\label{sec:degenerate}

In this section, with Lemma~\ref{lem:degdecay} we will establish an excess improvement estimate for the degenerate case which is characterized
by the fact that $\Phi(x_0,\varrho)$ is ``large'' compared to $\varphi(|(D{\bf u})_{x_0,\varrho}|)$. 

In view of Lemma~\ref{lem:lemma4.3}, this will be achieved via the $\varphi$-harmonic approximation lemma (Lemma~\ref{lem:phiharmapprox}) which allows to approximate the original minimizer by a $\varphi$-harmonic function. In this way, 
one can transfer the a priori estimates for $\varphi$-harmonic functions (Proposition~\ref{generalgrowth}) to the minimizer.


\begin{lemma}
Let $\gamma_0>0$ be the exponent of Proposition~\ref{generalgrowth}. Then, for every $0<\gamma<\gamma_0$ and every $\kappa,\mu\in(0,1)$ there exist $\varepsilon_\#, \tau\in(0,1)$ and $\varrho_\#\in(0,1]$ depending on $n,N,\mu_1,\mu_2, c_0, \beta_0$, $L,\nu,\gamma,\gamma_0,\mu$ and $\kappa$ ($\varepsilon_\#$ also depends on $\tau$ and $\sigma(\delta)$, where $\delta$ satisfies \eqref{eq:condelta} below, and $\varrho_\#$ also depends on $\omega$ and $\mathcal{V}$) with the following property: if
\begin{align}
\kappa & \varphi(|(D{\bf u})_{x_0,\varrho}|) \leq \Phi(x_0,\varrho)\leq \varepsilon_\# \label{(4.23ab)} 
\end{align}
for $B_\varrho(x_0)\subseteq\Omega$ with $\varrho\in(0,\varrho_\#]$, then
\begin{equation}
\Phi(x_0,\tau\varrho) \leq \tau^{2\gamma} \Phi(x_0,\varrho)\quad \mbox{ and }\quad \Theta(x_0,\tau\varrho)<\mu\,.
\label{(4.23d)}
\end{equation}
\label{lem:degdecay}
\end{lemma}

\proof
Without loss of generality, we assume that $x_0=0$ and, correspondingly, we use the abbreviations $\Phi(\varrho)=\Phi(0,\varrho)$, $\Psi(\varrho)=\Psi(0,\varrho)$, $\Theta(\varrho)=\Theta(0,\varrho)$. Let $0<\gamma<\gamma_0$ be fixed, $\tau\in(0,\frac{1}{2^{\mu_2}}]$ to be specified later, and we set $\varepsilon:=\tau^{2\gamma_0+n}$. Furthermore,
let $\delta_0=\delta_0(n,N,\varphi,\nu,L,\varepsilon)\in(0,1]$ be the constant according to the $\varphi$-harmonic
approximation (Lemma~\ref{lem:phiharmapprox}) with $\theta$ the exponent of higher integrability as in \eqref{(34Stroffo)}.

Now, we have to check that ${\bf u}$ complies with \eqref{(23Stroffo)}, in order to apply Lemma~\ref{lem:phiharmapprox}. From the Poincar\'e inequality, the shift change formula \eqref{(5.4diekreu)} with $\eta=\kappa$ and \eqref{(4.23ab)} we have
\begin{equation}
\Psi(\varrho)\leq c_{P}\dashint_{B_\varrho} \varphi(|D{\bf u}|)\,\mathrm{d}x\leq 2^{\mu_2-1}c_{P} (c_{\kappa}+1+{\kappa}^{-1})\Phi(\varrho)\leq c_{P}c(\kappa,\mu_2)\epsilon_\#\,,
\label{stimaPsi}
\end{equation}
where $c(\kappa,\mu_2):=2^{\mu_2-1}(c_{\kappa}+1+{\kappa}^{-1})>1$. An analogous computation and the concavity of $\varphi^{-1}$ give
\begin{equation*}
\begin{split}
E(B_\varrho):=\dashint_{B_\varrho} \varphi(|D{\bf u}|)\,\mathrm{d}x &\leq \varphi\left(\varphi^{-1}\left(\frac{c(\kappa,\mu_2)\varepsilon_\#}{\varrho}\right)\right) \\
&\leq \varphi\left(\frac{(c(\kappa,\mu_2))^{\frac{1}{\mu_1}}\varphi^{-1}(\varepsilon_\#)}{\varrho}\right)\,,
\end{split}
\end{equation*}
whence
\begin{equation}
\Theta(\varrho)\leq (c(\kappa,\mu_2))^{\frac{1}{\mu_1}}\varphi^{-1}(\varepsilon_\#)=: \tilde{c}(\kappa,\mu_1,\mu_2)\varphi^{-1}(\varepsilon_\#)\,.
\label{stimaTheta}
\end{equation}
In applying Lemma~\ref{lem:lemma4.3} we choose $\delta>0$ (which, in turn, determines $\sigma=\sigma(\delta)>0$ such that \eqref{(4.20)} holds) in such a way that
\begin{equation}
c_*\delta\leq \frac{\delta_0}{2}\,,
\label{eq:condelta}
\end{equation}
where $c_*$ is the constant of Lemma~\ref{lem:lemma4.3}. Then, we choose $\varepsilon_\#<1$ such that
\begin{equation*}
\varepsilon_\#\leq \min\left\{\frac{1}{c_{P}c(\kappa,\mu_2)}\varphi\left(\frac{\delta_0\sigma c_*}{4}\right), \varphi\left(\frac{\mu}{\tilde{c}(\kappa,\mu_1,\mu_2)}\right), \frac{\tau^n}{2 c_{\frac{1}{2}}}\varphi(\mu)\right\}\,,
\end{equation*}
so that, with \eqref{stimaPsi}-\eqref{stimaTheta}, we have
\begin{equation*}
\frac{\varphi^{-1}(\Psi(\varrho))}{\sigma}\leq\frac{\delta_0}{4} \quad \mbox{ and }\quad\Theta(x_0,\varrho)\leq\mu\,.
\end{equation*}
We also determine a radius $\varrho_\#\in(0,1]$ according to
\begin{equation*}
\left[\omega(\mu)^{1-\frac{1}{s}}+\mathcal{V}(\varrho_\#)^{1-\frac{1}{s}}\right]^{\beta_2}\leq\frac{\delta_0}{4}\,,
\end{equation*}
where $\beta_2$ is defined in Lemma~\ref{lem:lemma4.3}. 
Recalling the definition of $\widetilde{H}(\varrho)$, for $\varrho\leq \varrho_\#$ we then have
\begin{equation*}
[\widetilde{H}(\varrho)]^{\beta_2}\leq\frac{\delta_0}{4}\,.
\end{equation*}
For such choice of $\varepsilon_\#$ and $\varrho_\#$ it holds that
\begin{equation*}
c_*\left(\delta+[\widetilde{H}(\varrho)]^{\beta_2}+\frac{\varphi^{-1}(\Psi(\varrho))}{\sigma}\right) \leq \delta_0\,,
\end{equation*}
so that, by Lemma~\ref{lem:phiharmapprox}, there exists a unique $\varphi$-harmonic ${\bf w}\in W^{1,\varphi}(B_{\varrho/2},\R^N)$ with ${\bf w}={\bf u}$ on $\partial B_{\varrho/2}$ that satisfies
\begin{equation*}
\left(\dashint_{B_{\varrho/2}}|{\bf V}(D{\bf u})-{\bf V}(D{\bf w})|^{2\theta}\,\mathrm{d}x\right)^{\frac{1}{\theta}}\leq \tau^{2\gamma_0+n} \dashint_{B_{\varrho}}\varphi(|D{\bf u}|)\,\mathrm{d}x\,,
\end{equation*}
and
\begin{equation*}
\dashint_{B_{\varrho/2}}\varphi(|D{\bf w}|)\,\mathrm{d}x\leq c \dashint_{B_{\varrho/2}}\varphi(|D{\bf u}|)\,\mathrm{d}x\,.
\end{equation*}
Taking into account the higher integrability result of Lemma~\ref{lem:higint}, Lemma~\ref{corollary2.10} implies that
\begin{equation}
\dashint_{B_{\varrho/2}} |{\bf V}(D{\bf u})-{\bf V}(D{\bf w})|^{2}\,\mathrm{d}x\leq \tau^{2\gamma_0+n} \dashint_{B_{\varrho}}\varphi(|D{\bf u}|)\,\mathrm{d}x\,,
\label{(4.25celok)}
\end{equation}
and since $\tau<\frac{1}{2}$, from Proposition~\ref{generalgrowth} we also have
\begin{equation}
\mI{B_{\tau\varrho}} \kabs{\bfV(D{\bf w})-(\bfV(D {\bf w}))_{\tau\varrho}}^2 \, \mathrm{d}x \, \leq \, c \, \tau^{2 \gamma_0} \, 	\mI{B_{\varrho/2}} \kabs{\bfV(D{\bf w})-(\bfV(D {\bf w}))_{\varrho/2}}^2 \, \mathrm{d}x\,.
\label{(4.26celok)}
\end{equation}
Thus, with \eqref{(4.25celok)}-\eqref{(4.26celok)} we infer 
\begin{equation*}
\begin{split}
\Phi(\tau\varrho) & \leq 4 \mI{B_{\tau\varrho}} \kabs{\bfV(D{\bf u})-(\bfV(D {\bf w}))_{\tau\varrho}}^2 \, \mathrm{d}x \\
 & \leq 8 \mI{B_{\tau\varrho}} \kabs{\bfV(D{\bf u})-\bfV(D {\bf w})}^2 \, \mathrm{d}x + 8 \mI{B_{\tau\varrho}} \kabs{\bfV(D{\bf w})-(\bfV(D {\bf w}))_{\tau\varrho}}^2 \, \mathrm{d}x \\
& \leq c\tau^{-n}(\tau^{2\gamma_0+n}) \dashint_{B_\varrho}\varphi(|D{\bf u}|)\,\mathrm{d}x + c \tau^{2 \gamma_0} \, 	\mI{B_{\varrho/2}} \kabs{\bfV(D{\bf w})-(\bfV(D {\bf w}))_{\varrho/2}}^2 \, \mathrm{d}x \\
& \leq \tilde{c}_1\tau^{2\gamma_0}\dashint_{B_\varrho}\varphi(|D{\bf u}|)\,\mathrm{d}x
\end{split}
\end{equation*}
for some constant $\tilde{c}_1=\tilde{c}_1(n,N,\mu_1,\mu_2,c_1)>0$. Now, by virtue of the computation in \eqref{stimaPsi} we conclude that
\begin{equation*}
\Phi(\tau\varrho) \leq  \tilde{c}_1c(\kappa,\mu_2)\tau^{2\gamma_0} \Phi(\varrho)\,,
\end{equation*}
whence \eqref{(4.23d)} follows if we choose $\tau$ such that 
\begin{equation*}
\tau\leq \left(\frac{1}{\tilde{c}_1c(\kappa,\mu_2)}\right)^{\frac{1}{2(\gamma_0-\gamma)}}\,.
\end{equation*}
As for the second assertion in \eqref{(4.23d)}, the shift change formula \eqref{(5.4diekreu)} for $\eta=\frac{1}{2}$, with \eqref{(4.23ab)}, \eqref{stimaTheta} and the choice of $\varepsilon_\#$ shows that
\begin{equation*}
\begin{split}
E(B_{\tau\varrho}):=\dashint_{B_\tau\varrho} \varphi(|D{\bf u}|)\,\mathrm{d}x & \leq 2^{\mu_2-1}\left(c_{\frac{1}{2}}\tau^{-n}\dashint_{B_{\varrho}} \varphi_{|(D\bfu)_{\varrho}|}(|D{\bf u}-(D\bfu)_{\varrho}|)\,\mathrm{d}x + \frac{3}{2}\varphi(|(D\bfu)_{\varrho}|)\right) \\
&\leq 2^{\mu_2-1}\left(c_\frac{1}{2}\tau^{-n}\Phi(\varrho) +\frac{3}{2}E(B_{\varrho})\right)\\
& \leq 2^{\mu_2-1}\left(c_\frac{1}{2} \tau^{-n}\varepsilon_\#+\frac{3}{2}E(B_{\varrho})\right)\\
&  \leq 2^{\mu_2}\tau\varphi\left(\frac{\mu}{\tau\varrho}\right)
\end{split}
\end{equation*}
whence the assertion follows since $\tau\leq \frac{1}{2^{\mu_2}}$. 
\endproof

\subsection{Proof of Theorem~\ref{theorem-result-1}}

\begin{lemma}
Under the assumptions of Theorem~\ref{theorem-result-1}, let $\alpha\in(0,1)$. Then there exist constants $\varepsilon_\#$, $\delta_*$ and $\tilde{\varrho}$ such that the conditions
\begin{equation}
\Phi(x_0,\varrho)<\varepsilon_\# \quad \mbox{ and } \quad \Theta(x_0,\varrho)<\delta_*\,,
\label{(3.65verena)}
\end{equation}
for $B_\varrho(x_0)\subseteq\Omega$ with $\varrho\in(0,\tilde{\varrho}]$ imply
\begin{equation}
\bfu \in C^{0,\alpha}(\overline{B_{\varrho/2}(x_0)})\,.
\label{(3.66verena)}
\end{equation} 
\label{lem:lemma3.15}
\end{lemma}

\proof
Without loss of generality, we assume that $x_0=0$ and, correspondingly, we omit the dependence on it.
Let $\varepsilon_*,\delta_*$, $\vartheta\in(0,1)$ and $\varrho_*\in(0,1]$ be the constants of Lemma~\ref{lem:lemma3.13}. We then choose $\mu=\delta_*$ in Lemma~\ref{lem:degdecay} leaving $\kappa$ unchanged. This fixes the constants $\varepsilon_\#, \tau$ and $\varrho_\#$. We set $\tilde{\varrho}:=\min\{\varrho_*,\varrho_\#\}$.

We introduce the set of integers 
\begin{equation*}
\mathbb{S}:=\left\{k\in\mathbb{N}_0:\,\, \kappa \varphi(|(D{\bf u})_{\varrho}|) \leq \Phi(\tau^k\varrho)\right\}\,,
\end{equation*}
and we distinguish between the cases $\mathbb{S}=\mathbb{N}_0$ and $\mathbb{S}\neq\mathbb{N}_0$.\\
\noindent
\emph{The case $\mathbb{S}=\mathbb{N}_0$.} We prove by induction that the bounds
\begin{equation}
\Phi(\tau^k\varrho)<\varepsilon_\# \quad \mbox{ and }\quad \Theta(\tau^k\varrho)<\delta_*
\label{(Dk)}
\end{equation}
hold for every $k\in\N_0$. The case $k=0$ is trivial from the assumption \eqref{(3.65verena)}. Now, since $k\in\mathbb{S}=\N_0$, the assumption \eqref{(4.23ab)} of Lemma~\ref{lem:degdecay} hold with $\tau^k\varrho$ in place of $\varrho$. Then, an application of Lemma~\ref{lem:degdecay} gives \eqref{(Dk)} for $k+1$ (recall that $\tau<1$). The validity of \eqref{(Dk)} implies, as in Lemma~\ref{lem:lemma3.13}, that the Morrey-type estimate
\begin{equation}
\Theta(y,r)\leq c\delta_*\left(\frac{r}{\varrho}\right)^\alpha
\label{(5.10bis)}
\end{equation}
holds for every $\alpha\in(0,1)$, for all $y\in B_{\varrho/2}(x_0)$ and $r\in(0,\varrho/2]$. For $y,z\in B_{\varrho/2}$, with $|y-z|\leq\varrho/4$ we estimate the telescopic sum
\begin{equation*}
\frac{|\bfu(y)-\bfu(z)|}{|y-z|}\leq \sum_{j\in\mathbb{Z}}\frac{1}{2^j}\dashint_{B_{r_j}}\frac{|\bfu(x)-(\bfu)_j|}{|y-z|}\,\mathrm{d}x \leq  \sum_{j\in\mathbb{Z}}\frac{1}{2^j}\dashint_{B_{r_j}}\frac{|\bfu(x)-(\bfu)_j|}{r_j}\,\mathrm{d}x\,,
\end{equation*}
where $B_{r_j}:=B_{2^{1-j}|y-z|(y)}$ for $j\geq0$ and $B_{r_j}:=B_{2^{1+j}|y-z|(z)}$ for $j<0$. Now, with the Poincar\'e inequality we get
\begin{equation*}
\begin{split}
\frac{|\bfu(y)-\bfu(z)|}{|y-z|}\leq c\sum_{j\in\mathbb{Z}}\dashint_{B_{r_j}}|D\bfu|\,\mathrm{d}x &\leq c \sum_{j\in\mathbb{Z}}\varphi^{-1}\left(\dashint_{B_{r_j}}\varphi(|D\bfu|)\,\mathrm{d}x\right)\\
& \leq c \sum_{j\geq0} \frac{\Theta(y,2^{1-j}|y-z|)}{2^{1-j}|y-z|} + c \sum_{j<0} \frac{\Theta(z,2^{1+j}|y-z|)}{2^{1+j}|y-z|}\,.
\end{split}
\end{equation*}
Finally, with the estimate \eqref{(5.10bis)} we infer
\begin{equation}
\begin{split}
|\bfu(y)-\bfu(z)|&\leq c \delta_* |y-z|^\alpha\varrho^{-\alpha} \left(\sum_{j\geq0}2^{(\alpha-1)(1-j)}+\sum_{j<0}2^{(\alpha-1)(1+j)}\right)\\
& \leq c \delta_* |y-z|^\alpha\varrho^{-\alpha}\,,
\end{split}
\label{eq:holderianity}
\end{equation} 
which, in particular, implies that $\bfu\in C^{0,\alpha}(\overline{B_{\varrho/2}})$.
\\
\noindent
\emph{The case $\mathbb{S}\neq\mathbb{N}_0$.} In this case, there exists $k_0:=\min \mathbb{N}\backslash\mathbb{S}$. Since $k\in\mathbb{S}$ for any integer $k<k_0$ we can iterate as in the case $\mathbb{S}=\mathbb{N}_0$ for $k=0,1,\dots,k_0-1$ to infer that \eqref{(Dk)} holds for any $k\leq k_0$. By the definition of $\mathbb{S}$ we have
\begin{equation*}
\Phi(\tau^{k_0}\varrho)<\kappa\varphi(|(D{\bf u})_{\varrho}|)\,,
\end{equation*}
which together with the second inequality in \eqref{(Dk)} with $k=k_0$ ensures that the assumptions \eqref{eq:0step} of Lemma~\ref{lem:lemma3.13} are satisfied for $\varrho$ replaced by $\tau^{k_0}\varrho$. The, by virtue of this lemma, we have
\begin{equation}
\frac{\Phi(\vartheta^m\tau^{k_0}\varrho)}{{\varphi(|(D{\bf u})_{\vartheta^m\tau^{k_0}\varrho}|)}}\leq \kappa \quad \mbox{ and }\quad \Theta(\vartheta^m\tau^{k_0}\varrho)\leq\delta_*
\label{(Nl)}
\end{equation}
for every $m\in\mathbb{N}_0$.

Now, we consider an arbitrary radius $r\in(0,\varrho]$. If $r\in(\tau^{k_0}\varrho/2, \varrho]$ we find $0\leq k \leq k_0$ such that $\tau^{k+1}\leq r \leq \theta^k$ and then we can argue as in the case $\mathbb{S}=\mathbb{N}_0$. In the case $r\in(0,\tau^{k_0}\varrho/2]$, instead, we find $m\in\mathbb{N}_0$ such that $\vartheta^{m+1}\tau^{k_0}\varrho<r\leq \vartheta^{m}\tau^{k_0}\varrho$. Then arguing as in the proof of \eqref{(5.10Stroffo)} and taking into account the second estimate in \eqref{(Nl)}, we have

\begin{equation*}
r^{1-\alpha}\varphi^{-1}(E(B_r(y)))\leq c (\vartheta^m\tau^{k_0}\varrho)^{1-\alpha} \varphi^{-1}(E(B_{\vartheta^m\tau^{k_0}\varrho})) \leq \frac{c\delta_*}{(\vartheta^m\tau^{k_0}\varrho)^{\alpha}}
\end{equation*}
for every $y\in B_{\vartheta^m\tau^{k_0}\varrho/2}\subseteq B_{\varrho/2}$ whence
\begin{equation*}
\Theta(y,r)\leq c\frac{\delta_*}{(\vartheta^m\tau^{k_0})^{\alpha}}\left(\frac{r}{\varrho}\right)^\alpha\,.
\end{equation*}
At this point, we can argue as in the case $\mathbb{S}=\N_0$ for the proof of \eqref{eq:holderianity}, whence \eqref{(3.66verena)} follows thus concluding the proof.
\endproof

\proof[Proof of Theorem~\ref{theorem-result-1}] Let $\varepsilon_\#$, $\delta_*$ and $\tilde{\varrho}$ be the constants of Lemma~\ref{lem:lemma3.15}. Let $\Sigma_1$ and $\Sigma_2$ be defined as in the statement of Theorem~\ref{theorem-result-1}. Note that, by Lebesgue's differentiation theorem, it holds that $|\Sigma_1\cup\Sigma_2|=0$. Thus, we are reduced to show that each $x_0\in\Omega\backslash(\Sigma_1\cup\Sigma_2)$ belongs to the set
\begin{equation*}
\Omega_0:=\left\{z_0\in\Omega:\,\, \bfu\in C^{0,\alpha}(U_{z_0},\R^N)\,\,\mbox{ for every }\alpha\in(0,1)\,\,\mbox{ and for some }U_{z_0}\subset\Omega\right\}\,,
\end{equation*}
where $U_{z_0}$ is an open neighborhood of $z_0$. For this, let $x_0\in\Omega$ be such that both the conditions
\begin{equation}
\mathop{\lim\inf}_{\varrho\searrow 0}\dashint_{B_\varrho(x_0)}|{\bf V}_{|(D\bfu)_{x_0,\varrho}|}(D\bfu-(D\bfu)_{x_0,\varrho})|^2\,\mathrm{d}x=0\,\,\mbox{ and }\,\, m_{x_0}:=\mathop{\lim\sup}_{\varrho\searrow 0}|(D\bfu)_{x_0,\varrho}|<+\infty
\label{regularcond}
\end{equation}
hold. 

We set 
\begin{equation}
\sigma:=\min\left\{\frac{1}{c_\varphi c_\frac{1}{2}2^{\mu_2}}\varphi\left({\delta_*}\right),\frac{\varepsilon_\#}{c_\varphi}\right\}\,,
\label{eq:sigmacond}
\end{equation}
where the constants $c_\varphi$, $c_\frac{1}{2}$ are specified later and, correspondingly, with \eqref{regularcond} we can find a radius $\bar{\varrho}$ such that
\begin{equation}
\bar{\varrho}\leq \frac{\varphi(\delta_*)}{3\cdot 2^{\mu_2+1} \varphi(m_{x_0}+1)}
\label{eq:radius}
\end{equation}
and
\begin{equation}
\dashint_{B_{\bar{\varrho}}(x_0)}|{\bf V}_{|(D\bfu)_{x_0,\bar{\varrho}}|}(D\bfu-(D\bfu)_{x_0,\bar{\varrho}})|^2\,\mathrm{d}x\leq\sigma\quad \mbox{ and }\quad |(D\bfu)_{x_0,\bar{\varrho}}|\leq m_{x_0}+1\,.
\label{regcond2}
\end{equation}
Now, recalling that
\begin{equation*}
\Phi(x_0,\bar{\varrho})\leq c_\varphi \dashint_{B_{\bar{\varrho}}(x_0)}|{\bf V}_{|(D\bfu)_{x_0,\bar{\varrho}}|}(D\bfu-(D\bfu)_{x_0,\bar{\varrho}})|^2\,\mathrm{d}x
\end{equation*}
and observing that, as a consequence of the shift-change formula \eqref{(5.4diekreu)} with $\displaystyle\eta=\frac{1}{2}$, conditions \eqref{regcond2} imply
\begin{equation*}
\begin{split}
\dashint_{B_{\bar{\varrho}}(x_0)}\varphi(|D\bfu|)\,\mathrm{d}x& \leq 2^{\mu_2-1}\left(c_\frac{1}{2} c_\varphi \dashint_{B_{\bar{\varrho}}(x_0)}|{\bf V}_{|(D\bfu)_{x_0,\bar{\varrho}}|}(D\bfu-(D\bfu)_{x_0,\bar{\varrho}})|^2\,\mathrm{d}x+\frac{3}{2}\varphi(|(D\bfu)_{x_0,\bar{\varrho}}|)\right)\\
& \leq 2^{\mu_2-1}\left(c_{\frac{1}{2}}c_\varphi\sigma + \frac{3}{2}\varphi(m_{x_0}+1)\right) \leq \frac{1}{2}\varphi(\delta_*) + \frac{1}{2\bar{\varrho}}\varphi(\delta_*)\leq \varphi\left(\frac{\delta_*}{\bar{\varrho}}\right)\,,
\end{split}
\end{equation*}
with the choice of $\sigma$ \eqref{eq:sigmacond}, corresponding to the radius $\bar{\varrho}\in(0,\tilde{\varrho}]$ as in \eqref{eq:radius} there holds
\begin{equation*}
\Phi(x_0,\bar{\varrho})<\varepsilon_\#\quad \mbox{ and }\quad \Theta(x_0,\bar{\varrho})<\delta_*\,.
\end{equation*}
By the absolute continuity of the integral, we can find an open neighborhood $U_{x_0}$ of $x_0$ such that
\begin{equation*}
\Phi(x,\bar{\varrho})<\varepsilon_\#\quad \mbox{ and }\quad \Theta(x,\bar{\varrho})<\delta_*
\end{equation*}
for every $x\in U_{x_0}$. We can apply Lemma~\ref{lem:lemma3.15} at each point of $U_{x_0}$, proving that $\bfu\in C^{0,\alpha}(U_{x_0},\R^N)$ for every $\alpha\in(0,1)$. Thus, $x_0\in\Omega_0$ and the proof is concluded.
\endproof

\begin{acknow} 
G. Scilla has been supported by the Italian Ministry of Education, University and Research through the Project Variational methods for stationary and evolution problems with singularities and interfaces (PRIN 2017). The research of B. Stroffolini was supported by PRIN Project 2017TEXA3H.
\end{acknow}

\Addresses

\end{document}